\documentclass[12pt]{amsart}

\usepackage{amsmath}
\usepackage{latexsym}
\usepackage{amssymb}
\usepackage{amsfonts}
\usepackage[usenames]{color}
\usepackage{graphics, hyperref}

\renewcommand{\proof}{\par\noindent{\it Proof.\ \ }}
\def\qed{\ifmmode\square\else\nolinebreak\hfill
$\square$\fi\par\vskip12pt}

\def\ov{\overline} 
\def\l{\langle} \def\r{\rangle} 
\def\char{{\rm char}}

\def\calO{{\mathcal O}}
\def\ZZ{{\rm C}}

\def\div{\,\big|\,} \def\notdiv{{\,\not\big|\,}}

\def\FF{\mathbb F} \def\ZZ{\mathbb Z}

\def\PG{{\rm PG}} \def\Cay{{\sf Cay}} \def\Cos{{\sf Cos}}

\def\D{{\rm D}} \def\S{{\rm S}} 
 \def\M{{\sf M}} \def\soc{{\sf soc}}
 \def\C{{\bf C}} \def\N{{\bf N}}
\def\Z{{\bf Z}}  \def\mod{{\sf mod~}}
 \def\char{{\sf \,char\,}}
  
 \def\Aut{{\sf Aut}} \def\Inn{{\sf Inn}}

  \def\K{{\sf K}}
\def\Alt{{\sf Alt}}

\def\Ga{{\it \Gamma}}

\def\a{\alpha}   \def\s{\sigma}
\def\t{\tau} \def\o{\omega}

\def\GF{{\rm GF}}

 \def\PSp{{\rm PSp}}
 
\def\GammaL{{\rm \Gamma L}}

\def\A{{\rm A}}

\def\PSL{{\rm PSL}}\def\PGL{{\rm PGL}}
\def\GL{{\rm GL}}

\def\AGL{{\rm AGL}}

 \def\PSU{{\rm PSU}}

  \def\D{{\rm D}}

\newtheorem{theorem}{Theorem}[section]
\newtheorem{lemma}[theorem]{Lemma}%
\newtheorem{example}[theorem]{Example}%
\newtheorem{construction}[theorem]{Construction}%

\def\FF{\mathbb F}    \def\ZZ{\mathbb Z}

\def\mod{{\sf mod~}}

\def\K{{\bf K}} 
 
\def\C{{\bf C}}\def\N{{\bf N}}

\def\Cay{{\sf Cay}} \def\Cos{{\sf Cos}}
\def\Aut{{\sf Aut}} \def\Inn{{\sf Inn}}

\def\PG{{\rm PG}}
\def\D{{\rm D}} 
\def\S{{\rm S}} 
 \def\M{{\rm M}}

\def\soc{{\sf soc}}

\def\a{\alpha}   \def\s{\sigma}
\def\t{\tau}

\begin{document}

\title[Edge-transitive Cayley graphs]
{Tetravalent edge-transitive Cayley graphs of Frobenius groups}

%\thanks{1991 MR Subject Classification 20B15, 20B30}%, 05C25.}
%\thanks{This work forms part of the PhD project of the second-named author;
%partially supported by an NNSF(K1020261), and an ARC Discovery Grant.}

\author{Lei Wang}
\address{Beijing International Center for Mathematical Research\\
Peking University\\
Beijing, 100871\\
P. R. China}
\email{wanglei-2468@163.com(Lei Wang)}

\date\today

\maketitle

\begin{abstract}
In this paper, we give a characterization  for a class of edge-transitive Cayley graphs, and provide methods for
constructing Cayley graphs with certain symmetry properties.
Also this study leads to construct and characterise a new family of half-transitive graphs.
\end{abstract}

\qquad {\textsc k}{\scriptsize \textsc {eywords.}} {\footnotesize
Frobenius group, Edge-transitive graph, Coset graph, Cayley graph}

\section{Introduction}
Graphs considered in this paper are assumed to be finite, simple, and unless stated otherwise,
connected and undirected. For a graph $\Ga$, let $V\Ga$, $E\Ga$ and $\Aut\Ga$ denote its vertex set, edge set and
the full automorphism group, respectively.
 If there exists a subgroup $X\leqslant\Aut\Ga$ is transitive on $V\Ga$ or $E\Ga$,
then the graph $\Ga$ is said to be $X$-vertex transitive or $X$-edge transitive,
respectively. A sequence $v_0,v_1,\ldots, v_s$ of vertices of $\Ga$ is called an $s$-arc if
$v_{i-1}\not=v_{i+1}$ for $1\leqslant i\leqslant s-1$, and $\{v_i,v_{i+1}\}$ is an edge for $0\leqslant i\leqslant s-1$.
The graph $\Ga$ is called $(X, s)$-arc-transitive, if $X$ is transitive on the $s$-arcs of $\Ga$;
if in addition $X$ is not transitive on the $(s+1)$-arcs,
then $\Ga$ is said to be $(X,s)$-transitive. In particular,
a $1$-arc is simply called an arc, and an $(X,1)$-arc-transitive graph is called $X$-arc transitive.

A graph $\Ga$ is called a Cayley graph if there exist a group $G$ and
a subset $S\subset G\setminus \{1\}$ with $S=S^{-1}{:}=\{g^{-1}\mid g\in S\}$
such that the vertices of $\Ga$ may be identified with the elements of $G$
in such a way that $x$ is adjacent to $y$ if and only if $yx^{-1}\in S$.
The Cayley graph $\Ga$ is denoted by $\Cay(G, S)$.
Throughout this paper, denote by {\bf1} the vertex of $\Cay(G,S)$ corresponding
to the identity of $G$.

It is well-known that a graph $\Ga$ is a Cayley graph of a group $G$ if and only if the full automorphism
group $\Aut\Ga$ contains a subgroup which is regular on vertices and
isomorphic to $G$. In particular, a Cayley graph $\Cay(G,S)$ is vertex-transitive of order $|G|$.
However, a Cayley graph  is of course not necessarily edge-transitive.
Thus, characterizing the Cayley graphs which are edge-transitive is a current hot topic in algebra graph theory.
For instance, see \cite{Xin gui,Hua Zhang,Xiuyun1,Ming Yao} for those with valency $4$,
see \cite{Zai Ping Lu} for a classification of connected edge-transitive tetravalent Cayley graphs of
square-free order, and  \cite{Corr} for a classification of
normal edge-transitive Cayley graphs of Frobebius groups of order a product of two primes.
In this paper, a characterization is given of tetravalent edge-transitive Cayley graphs of
a class of primitive Frobenius groups.
This study provides a method for constructing edge-transitive graphs of valency 4,
and is then applied to construct a new family of half-transitive graphs.
To state this result, we need more definitions.

For an $X$-vertex-transitive graph $\Ga$ and a normal subgroup $N\lhd X$,
the normal quotient graph $\Ga_N$ induced by $N$
is the graph which has vertex set $V\Ga_N=\{u^N\div u\in V\Ga\}$ such that $u^N$
and $v^N$ are adjacent if and only if $u$ is adjacent in $\Ga$ to some vertex in $v^N$.
Furthermore, if the valency of $\Ga_N$ equals the valency of $\Ga$, then $\Ga$ is called a normal cover of $\Ga_N$.

For an integer $m\geqslant3$, we denote by $\C_{m[2]}$ the
lexicographic product of the empty graph $2\bf K_1$ of order $2$ by a cycle $\C_m$ of size $m$,
which has vertex set $\{(i,j)\div  1\leqslant i\leqslant m, 1\leqslant j\leqslant 2\}$ such that
$(i,j)$ and $(i',j')$ are adjacent if and only if $i-i'\equiv\pm1\,(\mod m)$.

A group $G$ is said to be a {\it Frobenius group} if and only if
$G$ has the form $G=W{:}H$ such that each non-identity element of $H$
centralises no non-identity element of $W$, that is, $xy\not=yx$
for any $x\in W\setminus\{1\}$ and $y\in H\setminus\{1\}$.
In particular, $G$ is called a {\it primitive Frobenius group} if $H$ acts irreducibly on $W$,
refer to \cite{DM-book}.

%For positive integers $a,m$, a divisor $r$ of $a^m-1$ is called a {\it primitive divisor} if
%$r$ does not divide $a^i-1$ for $i<m$.

Let $\FF$ be a field, $R$ be a group and $V$ be an $\FF R$-module.
Suppose that $V=V_1\oplus\cdots\oplus V_r$ $(r>1)$, where $V_i$ are subspaces of $V$ which are
transitively permuted by the action of $R$.
We call $R$ imprimitive on $V$ if there exists such decomposition.
Otherwise, $R$ is called  primitive on $V$.

\begin{theorem}\label{soluble}
Let $G=W{:}H\cong\ZZ_p^d{:}\ZZ_n$ be a primitive Frobenius group,
where $d,n$ are integers, and $p$ is a prime.
Assume that $\Ga$ is a connected tetravalent $X$-edge-transitive
Cayley graph of $G$, where $G\leqslant X\leqslant \Aut\Ga$.
If $X$ is soluble, then one of the following statements holds:
\begin{itemize}
\item[(1)]
$G$ is normal in $X$, and $X_1\leqslant\D_8$;
\item[(2)]
$G\cong\D_{2p}$,  $\Ga\cong\C_{p[2]}$, and $\Aut\Ga\cong\ZZ_2^p{:}\D_{2p}$;

\item[(3)]
$X=W{:}((N{:}H).O)$ with $\soc(X)=W\times L$, and $X_1=N.O$,
where  $N\cong\ZZ_2^l$ with $2\leqslant l\leqslant d$, $L\cong1$ or $\ZZ_2$, and $O\cong1$ or $\ZZ_2$,
satisfying the following statements:

\begin{itemize}
\item[(a)]
there exist $x_1,\ldots,x_d\in W$ and $\tau_1,\tau_2,\ldots,\tau_d\in N$
such that $W=\l x_1,x_2,\ldots,x_d\r$,
$\l x_i,\tau_i\r\cong\D_{2p}$ and $N=\l\tau_i\r\times\C_N(x_i)$ for $1\leqslant i\leqslant d$;

\item[(b)]
$H$ does not centralise $N$, and $H$ is imprimitive on $W$;

\item[(c)]
$X/(WN)\cong\ZZ_n$ or $\D_{2n}$, and $\Ga$ is $X$-arc-transitive if and only if
$X/(WN)\cong\D_{2n}$;

\end{itemize}

\item[(4)]
$\Ga_W\cong\C_{\frac{n}{2}[2]}$, $\Ga$ is a cover of $\Ga_W$  and $X=W{:}((NH).O)$ such that
\begin{itemize}
\item[(i)]
$X_1\leqslant N.O$, $N\cap H\cong\ZZ_2$, and $H$ normalizes $N$, but $H$ does not centralise $N$,
where $N\cong\ZZ_2^l$ with $2\leqslant l\leqslant \frac{n}{2}$, $4$ divides $n$, and $O\cong1$ or $\ZZ_2$;

\item[(ii)]
$W$ is the unique minimal normal subgroup of $X$, and $H$ is imprimitive on $W$;

\item[(iii)]
$X/(WN)\cong\ZZ_{\frac{n}{2}}$ or $\D_n$, and $\Ga$ is $X$-arc-transitive if and only if
$X/(WN)\cong\D_n$;

\end{itemize}

\item[(5)]
$X=((WN){:}H).O$
and $\Ga$ is $X$-arc-transitive if and only if $X/(WN)\cong\D_{2n}$,
where  $W\cong\ZZ_2^d$, $N$ is a $2$-group, and $O\cong1$ or $\ZZ_2$.
\end{itemize}
\end{theorem}
\vskip0.1in
{\noindent\bf Remarks on Theorem~\ref{soluble}.}
\begin{itemize}
\item[(a)]
The Cayley graph $\Ga$ in part~(1), called a normal edge-transitive graph, is studied in \cite{Praeger}.
Furthermore, if $X=\Aut\Ga$,  then $\Ga$ is called a normal Cayley graph, introduced in \cite{Xu}.

\item[(b)]
$H$ acts irreducibly on $W$ if and only if
$n$ does not divide $p^m-1$ for any proper divisor $m$ of $d$
(such $n$ is called a {\it primitive divisor} of $p^d-1$), refer to \cite[Proposition~2.3]{DF}.

\item[(c)]
Lemma~\ref{Lee} and Lemma~\ref{Le} show that every group $X$ satisfies part~(3) or part~(4)
if and only if $H$ is imprimitive on $W$, see Constructions~\ref{im} and~\ref{yin}.
In addition, $H$ is imprimitive on $W$ if and only if there exists some prime $k$ dividing $d$
such that $n$ divides $k(p^{\frac{d}{k}}-1)$, see \cite[Proposition~2.8]{DF}.
\end{itemize}
\vskip0.1in

\begin{theorem}\label{isoluble}
Using the notation defined in Theorem~$\ref{soluble}$,
if $X$ is insoluble, then one of the following holds:
\begin{itemize}

\item[(1)]
$G\cong\ZZ_p^4{:}\ZZ_5$, $X=W.\overline X$ and $\Ga_W\cong\K_5$, where $\soc(\ov X)\cong\A_5$, and
$\Ga$ is constructed as in Construction~$\ref{AA}$;

\item[(2)]
$G\cong\ZZ_p^4{:}\ZZ_{10}$, $X=W.(\overline X\times\ZZ_2)$, and $\Ga_W\cong\K_{5,5}-5\K_2$, where $\soc(\ov X)\cong\A_5$,
and $\Ga$ is constructed as in
Construction~$\ref{AAAAA}$;

\item[(3)]
$\Ga$ is isomorphic to one of the graphs listed in Table~$1$.
\end{itemize}
\end{theorem}
\vskip0.07in
\centerline{\bf TABLE~1:\ {\rm Graphs which are not normal edge-transitive.}}
%Table~1 Almost simple automorphism groups.
\[\begin{array}{|c|c|c|c|c|}\hline
\Aut\Ga& G & (\Aut\Ga)_1  & \Ga \\ \hline
\PSL(3,3){:}\ZZ_2 &\D_{26}  & \ZZ_3^2{:}\GL(2,3) & \mbox{Example}~\ref{TA}\\ \hline
\PGL(2,7) & \D_{14} & \S_4 & \mbox{Example}~\ref{TA}\\\hline
\PGL(2,7) & \ZZ_7{:}\ZZ_3 & \D_{16}&  \mbox{Example}~\ref{A} \\\hline
\PGL(2,7) & \ZZ_7{:}\ZZ_6 & \D_8  & \mbox{Example}~\ref{B}\\\hline
%\PSL(2,11) &\ZZ_{11}{:}\ZZ_{5}  & \A_4 & \mbox{Example}~\ref{C}\\\hline
\PSL(2,23) &\ZZ_{23}{:}\ZZ_{11}  & \S_4 & \mbox{Example}~\ref{C}\\ \hline
\PGL(2,11) &\ZZ_{11}{:}\ZZ_{5}  & \S_4 &  \mbox{Example}~\ref{C}\\ \hline
\PGL(2,11)\times\ZZ_2 &\ZZ_{11}{:}\ZZ_{10}  & \S_4 & \mbox{Example}~\ref{D}\\\hline
\end{array}\]

\vskip0.2in
A graph is said to be half-transitive if its automorphism group acts transitively on the vertex set and
edge set but intransitively on the arc set.
Constructing and characterising half-transitive graphs
was initiated by Tutte (1965), and is a currently active topic, see \cite{Maru,D.Maru,D.Mar,DD} for references.
Theorem~\ref{soluble} provides a method for characterising some classes
of half-transitive graphs of valency $4$. The following theorem is such an example.

\begin{theorem}\label{solubles}
Let $G=W{:}\langle h\rangle\cong\ZZ_p^d{:}\ZZ_n$ be a primitive Frobenius group,
where $d>1$ is odd, $p$ is an odd prime, and
$n$ is an integer.
Let $\Ga$ be a connected tetravalent edge-transitive
Cayley graph of $G$. Assume that $\l h\r$ is primitive on $W$.
Then $\Aut\Ga=G{:}\ZZ_2$, $\Ga$ is half-transitive, and
$\Ga\cong\Ga_i=\Cay(G,S_i)$, where $1\leqslant i\leqslant\lfloor\frac{n-1}{2}\rfloor$,
$(n,i)=1$, and
\[S_i=\{ah^i, a^{-1}h^i, (ah^i)^{-1}, (a^{-1}h^i)^{-1}\}, where\,\, a\in W\setminus \{1\}.\]
Moreover, if $p^ri\equiv\pm j\ (\mod n)$ for some $r\geqslant 0$, then $\Ga_i\cong\Ga_j$.
\end{theorem}

\section{Preliminary results}
In this section, we quote some preliminary results, which will be used in the subsequent sections.

\vskip0.07in
For a core-free subgroup $H$ of $X$ and an element $g\in X\setminus H$, let
$[X{:}H]:=\{Hx\mid x\in X\}$, and define the coset graph
\[\Ga=\Cos(X,H, H\{g,g^{-1}\}H)\]
with vertex set $[X{:}H]$ such that $Hx$ and $Hy$ are adjacent
whenever $yx^{-1}\in H\{g,g^{-1}\}H$.  Then $\Ga$ is well-defined,
and $X$ induces a subgroup of $\Aut\Ga$ acting on $[X{:}H]$ by right multiplication,
namely, $\a : Hx\to Hxa$ for $x,a\in X$.
Label $v,w$ the two vertices of $\Ga$ corresponding to $H$ and $Hg$, respectively. Then
we have the following lemma.

\begin{lemma}\label{Cos}

For a coset graph
$\Ga=\Cos(X,H, H\{g,g^{-1}\}H)$, we have
\begin{itemize}
\item[(a)]
$\Ga(v)=\{Hgh|h\in H\}\cup\{Hg^{-1}h|h\in H\}$;
\item[(b)]
$\Ga$ is $X$-edge-transitive and $X$ is transitive on the vertices of $\Ga$;
\item[(c)]
$\Ga$ is connected if and only if $X=\l H,g\r$;
\item[(d)]
$H\cap H^g=X_{vw},$ the stabilizer of the arc $(v,w)$, where $H^g$
is the conjugate of $H$ by $g$;
\item[(e)]
the valency of $\Ga$ equals

\begin{equation*}
{val(\Ga)=}
\begin{cases}
|H{:}H\cap H^g| & \text{if $HgH=Hg^{-1}H$,}\\
2|H{:}H\cap H^g| &\text{otherwise};
\end{cases}
\end{equation*}

\item[(f)]
$\Ga$ is $X$-arc-transitive if and only if $HgH=Hg^{-1}H,$
which yields that $HgH=HoH$ for some $(2$-element $)$ $o\in\N_{X}(X_{vw})\setminus H$
with $o^2\in X_{vw}$ (refer to \cite{Zai Ping Lu}). (An element o in the group {\bf $X$}
is a $2$-element if its order is a power of {\bf $2$}).

\end{itemize}

Moreover, for any $X$-edge-transitive graph $\Sigma$,
if $X$ is transitive on $V\Sigma$, then the map $u^x\to Hx$ with $x\in X$ gives an isomorphism
from $\Sigma$ to $\Cos(X, H, H\{g, g^{-1}\}H)$,
where $u\in V\Sigma$, $H = X_u$ and $g\in X \setminus H$ with $u^g\in\Ga(u)$.

\end{lemma}

The vertex stabilizer for $s$-arc-transitive graphs of valency $4$ is known (refer to \cite{Weiss}).

\begin{lemma}\label{Lu}
Let $\Ga=(V\Ga,E\Ga)$ be a connected $(X,s)$-transitive graph of valency $4$.
Then $s$ and the stabilizer $X_1$ are listed in the following
table,
\[\begin{array}{|c|c|c|c|c|c|c|}\hline
s & 2 & 3 & 4 & 7   \\ \hline
X_1  & \A_4,\,\S_4 & \ZZ_3\times\A_4, (\ZZ_3\times\A_4).\ZZ_2,\,\S_3\times\S_4  &\ZZ_3^2{:}\GL(2,3)& [3^5]{:}\GL(2,3) \\
\hline
\end{array}\]
where $[3^5]$ is a $3$-group of order $3^5.$
\end{lemma}

Let $\Ga=(V\Ga,E\Ga)$ be a connected graph. Assume that $X\leqslant\Aut\Ga$ is transitive on both
$V\Ga$ and $E\Ga$. % and let the  $\Ga_N$, then we have
Then we have an important conclusion in the next lemma.

\begin{lemma}\label{insoluble}
Let $N\lhd X$.
If $\Ga$ is of valency $4$ and $X/N$ is insoluble, then $\Ga$ is a normal $N$-cover of $\Ga_N$.
\end{lemma}
\proof Pick any vertex $u\in V\Ga$. Let $B$ be an orbit of $N$ acting on $V\Ga$, which contains $u$.
By Lemma~\ref{Lu}, the  stabilizer $X_u$ is a $\{2,3\}$-group. In particular, $X_u$  is soluble.
Let $K$  be the kernel of $X$ acting on $\Ga_N$.
Then $K_u\unlhd X_u$, so $K_u$ is soluble.
Since $N$ is transitive on $B$, we have $K=NK_u$.
Note that $K/N\cong NK_u/N\cong K_u/(N\cap K_u)$,
$K/N$ is soluble.
Then $X/K\cong (X/N)/(K/N)$ is insoluble because $X/N$ is  insoluble.
So $\Aut\Ga_N$ is also insoluble, hence $\Ga_N$ is not  a cycle.
Since $\Ga$ is connected and the valency of $\Ga_N$ is a divisor
of the valency of $\Ga$, we conclude that $\Ga_N$ is of valency $4$,
and the lemma holds.
\qed

For a normal edge-transitive Cayley graph $\Ga=\Cay(G,S)$, let $\Aut(G,S)=\{\a\in\Aut(G)\div S^\a=S\}$,
we have a simple lemma.

\begin{lemma}\label{G}
Let $G=W{:}\l h\r\cong\ZZ_p^d{:}\ZZ_n$ be a primitive Frobenius group,
where $d, n$ are integers, and $p$ is a prime.
Let $\Ga=\Cay(G,S)$ be connected of valency $4$.
Assume that $\Aut\Ga$ has a subgroup $X$ such that
$\Ga$ is $X$-edge-transitive and $G\unlhd X$.
Then $X_1\leqslant\D_8$.
\end{lemma}
\proof Since $\Ga$ is connected, we have $\l S\r=G$, and so $\Aut(G,S)$ acts faithfully on $S$.
Hence $\Aut(G,S)\leqslant\S_4$. By \cite[Lemma~2.1]{Godsil}, we obtain $X\leqslant\N_{\Aut\Ga}(G)=G{:}\Aut(G,S)$.
So $X_1\leqslant\Aut(G,S)\leqslant\S_4$. Suppose that $3$ divides $|X_1|$. Then $X_1$ is $2$-transitive on $S$.
Hence $\Ga$ is $(X,2)$-arc-transitive, and
all elements in $S$ are involutions, see for example~\cite{flag}.

Pick any $s\in S$. Write $s=\s h^i$ where $\s\in W$ and $i$ is an integer.
Recall that $s$ is an involution,  we obtain that $h^{2i}=1$.
By \cite[Proposition~12.10]{Doerk},  $\Aut(G)\cong\ZZ_p^d{:}\GammaL(1,p^d)$.
For a finite group $T$,
it is known that  the action of $\Aut(T)$ on $T/\Z(T)$ is permutationally
isomorphic to the conjugation action of $\Aut(T)$ on $\Inn(T)$.
Since $G\cong\Inn(G)$, it follows from the above fact that we may identify $G$ with $\Inn(G)$ a normal subgroup of $\Aut(G)$.
Then write $\Aut(G)=W{:}M.L$, where $M\cong\ZZ_{p^d-1}$, and $L\cong\ZZ_d$.
Without loss of generality, we may assume that $h$ belongs to $ M$, refer to \cite[Proposition~2.5]{DF}.
For that case, $(h^i)^\eta=h^i$ for any $\eta\in M.L$.
Take any $\theta\in\Aut(G)$, $\theta$ has the form $xyz$, where $x\in W$, $y\in M$, and $z\in L$.
By easy calculations,  we have $s^\theta=\ov s h^i$, where $\ov s\in W$.
It follows that for each $a\in S$, $a$ has the form $\ov a h^i$ with $\ov a\in W$ because $X_1$ is transitive on $S$.
Recall that $\l S\r=G$, we have $h=(\s_1h^i)(\s_2h^i)\cdots(\s_mh^i)$ where $\s_j\in W$ for each $j$.
Since $W\unlhd G$ and $h\not=1$, we obtain $h=h^i$.  Consequently, $\l h\r\cong\ZZ_2$. By the definition of $G$,
we have $G\cong\D_{2p}$, and thus $\Aut(G)\cong\ZZ_p{:}\ZZ_{p-1}$.
However, since $X_1$ is $2$-transitive on $S$, we conclude that $X_1\cong\A_4$ or $\S_4$,
which is impossible.
Therefore, $X_1\leqslant\D_8$.
\qed

Finally, we quote a result about simple groups, which will be used later.

\begin{lemma}{\rm ({Kazarin\,\cite{L.Kazarin})}}\label{Kazarin}
Let $T$ be a non-abelian simple group which has a $2'$-Hall
subgroup. Then $T=\PSL(2,p)$, where $p=2^e-1$ is a prime.
Furthermore, $T=GH,$ where $G=\ZZ_p{:}\ZZ_{\frac{p-1}{2}}$ and
$H=\D_{p+1}=\D_{2^e}$.
\end{lemma}

\section{existence of graphs satisfying Theorem~$\ref{soluble}$ and Theorem~$\ref{isoluble}$}

%Before stating the next results, we make a convention that Frobenius groups
%(appearing in Theorems~\ref{soluble}-\ref{isoluble}) are referred to as primitive Frobenius groups.
In this section, we first construct some examples of graphs satisfying Theorem~\ref{soluble}.

\vskip0.07in
The following construction produces normal edge-transitive Cayley graphs admitting
a group $X$ satisfying part~(1) of Theorem~\ref{soluble}.

\begin{construction}\label{BBB}
{\rm
Let $p\geqslant5$ be a prime such that $p$ is a primitive divisor of $2^{p-1}-1$.
Let $G=W{:}\l h\r\cong\ZZ_2^{p-1}{:}\ZZ_p$ be a Frobenius group.
By \cite[Proposition~12.10]{Doerk}, we have $\Aut(G)\cong\ZZ_2^{p-1}{:}\GL(1,2^{p-1}).\ZZ_{p-1}$,
where $\ZZ_{p-1}$ is the group of Frobenius automorphisms.
Arguing similarly as Lemma~\ref{G},
%Since $G\cong\Inn(G)$, we may identify $G$ with $\Inn(G)$ a normal subgroup of $\Aut(G)$.
we may write $\Aut(G)=W{:}M{:}L$, and $h$ belongs to $M$,
where $M\cong\ZZ_{2^{p-1}-1}$, and $L\cong\ZZ_{p-1}$.
%Without loss of generality, we may assume that $\l h\r$ is contained in $M$ (refer to \cite[Proposition~2.4]{DF}).
For this case,  we may
identify $W$ with a field $\FF{:}=\FF_{2^{p-1}}$ of order $2^{p-1}$ and there exists $\a\in\FF$ of order
$p$ such that $\l h\r$ acts on each $x\in W$ by $h : x=\a x$.
By the definition,  $G$ is a primitive Frobenius group.

Let $\FF^{\#}=\l\o\r$, and let
$\s$ be a Frobenius automorphism of order $2$.
Then $\o^\s=\o^{2^{\frac{p-1}{2}}}$.
Let $X=G{:}\l\s\r$, and let $g=\o^{2^{\frac{p-1}{2}}+1} h$.
Set \[\Ga(2,p-1,p)=\Cos(X,\l\s\r,\l\s\r\{g,g^{-1}\}\l\s\r).\]
}
\end{construction}

\begin{lemma}\label{BB}
Let $\Ga=\Ga(2,p-1,p)$ be a graph constructed in Construction~$\ref{BBB}$.
Then $\Ga$ is a connected normal $X$-edge-transitive Cayley graph of $G$ of valency $4$.
\end{lemma}
\proof By the definition, $\l\s\r$ is core-free in $X$, and hence $X\leqslant\Aut\Ga$.
Now $X=G\l\s\r$ and $G\cap\l\s\r=1$, and thus $G$ is regular on the vertex set $[X{:}\l\s\r]$.
So $\Ga$ is a Cayley graph of $G$, which has order $2^{p-1}p$.

Let $Y=\l g,\s\r$.
Noting that $p\geqslant 5$, we conclude that $2^{\frac{p-1}{2}}\not\equiv-1(\mod 2^{p-1}-1)$.
It implies that $\o^\s\not=\o^{-1}$.
However, since $\a^\s=\a^{-1}$, we have that
$h^\s=h^{-1}$. Furthermore, $\s$ induces an automorphism of $G$.
Then we have \[g^\s=(\o^{2^{\frac{p-1}{2}}+1} h)^\s=(\o^\s)^{2^{\frac{p-1}{2}}+1} h^\s=\o^{2^\frac{p-1}{2}+1}h^{-1}.\]
Let $\ov g=g^\s g$. Then $\ov g=\o^{2(2^{\frac{p-1}{2}}+1)}\a$.
Denote by $\ell$ the integer $2^{\frac{p-1}{2}}-1$.
Recall that $p$ is a primitive divisor of $2^{p-1}-1$,
we conclude that $(p,\ell)=1$. Thus $\ov g^\ell=\a^\ell$ belongs to $Y$. So does $\a$.
Consequently, $\o^{2(2^{\frac{p-1}{2}}+1)}$ belongs to $Y$, and so
$h$ belongs to $Y$. Since $\l h\r$ acts irreducibly on $W$,
we obtain that $X=Y$.
Thus $\Ga$ is connected.
It is  straightforward to show that $\l\s\r\cap\l\s\r^g=1$, and
hence $\l\s\r\cap\l\s\r^g$ has index $2$ in $\l\s\r$.
Since $X\leqslant\Aut\Ga$, it follows that $\Ga$ is not a cycle.
By Lemma~$\ref{Cos}$, $\Ga$ is connected, $X$-edge-transitive
and of valency $4$.
\qed

{\noindent\bf Remark.} In fact, the graphs in Construction $\ref{BBB}$ really exist.
For example, $p=5,11,13,19$, and so on.

\vskip0.07in

The following construction produces edge-transitive graphs admitting
a group $X$ satisfying part~(3) of Theorem~\ref{soluble} with $L\cong\ZZ_2$, and $O=1$.

\begin{construction}\label{im}
{\rm Let $X=W{:}(N{:}\l h\r)\cong\ZZ_p^d{:}(\ZZ_2^d{:}\ZZ_n)$, where
$p=2^\ell m+1$ be an odd prime,  $m\geqslant3$ is an odd number, and $\ell\geqslant1$,
such that $W\cong\ZZ_p^d$, $N\cong\ZZ_2^d$, and $\l h\r\cong\ZZ_n$ satisfy
\begin{itemize}
\item[(a)]
$d>1$,  $d$ divides $m$, and $2md$ is a primitive divisor of $p^d-1$;

\item[(b)]
$W=\prod_{i=1}^d\l x_i\r$, where $x_i=(1,\ldots,1,x,1,\ldots,1)$ with $o(x)=p$ for each $i$;

\item[(c)]
$N=\prod_{i=1}^d\l \t_i^{\frac{p-1}{2}}\r$, where $\t_i=(1,\ldots,1,\t,1,\ldots,1)$
with $x^{\t}=x^r$ and $r^{p-1}\equiv1\,(\mod p)$ for each $i$;

\item[(d)]
$h=c_1\t_1^{\frac{p-1}{2m}}(12\ldots d)$, where $c_1=(c,1,\ldots,1)$ with $x^c=x$, $\t^x=\t$, and $o(c)=2$.
\end{itemize}

Let $y=(x_1h)^{-1}$. Set
 \[\Ga(p,2,n)=\Cos(X,N,N\{ y,y^{-1}\}N).\]
}
\end{construction}

\begin{lemma}\label{Liu}
Let $\Ga=\Ga(p,2,n)$ be a graph constructed in Construction~$\ref{im}$, and let $G=W{:}\l h\r$. Then $\Ga$
is a connected tetravalent $X$-edge-transitive Cayley graph of Frobenius group $G$,
and $G$ is not  normal in $X$.
\end{lemma}
\proof By the definition, $N$ is core-free in $X$, and hence $X\leqslant\Aut\Ga$.
 Now $X=GN$ and $G\cap N=1$, and thus $G$ acts regularly on the vertex set $[X{:}N]$.
So $\Ga$ is a Cayley graph of $G$. Obviously, $G$ is not normal in $X$.

Let $H=\l h\r$. Suppose that $C{:}=\C_H(W)\not=1$. Then $C=\l h^\ell\r$, where $\ell$ divides $2md$.
Write $\ell=l_1d+l$, where $0\leqslant l_1<2m$, and $0\leqslant l<d$.
Let $\ov\t=\t_1\t_2\cdots\t_d$, and $\ov c=c_1\cdots c_d$, where $c_i=(1,\ldots,1,c,1,\ldots,1)$ for each $i$.
If $l=0$, then $h^\ell=\ov c~\ov \t^{\frac{p-1}{2m}l_1}$ and so
$x_1^{h^\ell}\not=x_1$, a contradiction.
Thus $l\not=0$.
Then $h^\ell=\ov c~\ov\t^{\frac{p-1}{2m}l_1}h^{l}$.
Let $(l,d)=k$.
Let $k'=l/k$, and $d'=d/k$. Relabeling if necessary, we may rewrite $\{1,\ldots,d\}=\{1_1,\ldots,i_j,\ldots,k_{d'}\}$.
Without loss of generality, we may assume that
$h^l=h_1\cdots h_k$,
where $h_i=(c_{i_1}\t_{i_1}^{\frac{p-1}{2m}})(c_{i_2}\t_{i_2}^{\frac{p-1}{2m}})\cdots(c_{i_{k'}}\t_{i_{k'}}^{\frac{p-1}{2m}})(i_1 i_2\cdots i_{d'})$.
Then $x_{1_1}^{h^\ell}=x_{1_2}^{r^{\frac{p-1}{2m}(l_1+1)}}\not=x_{1_1}$ because $1_1\not=1_2$, a contradiction.
Thus $H$ acts faithfully  on $W$.

We claim that $H$ is fixed-point-free on $W$.
Let $U=\l w\div w^h=w, w\in W\r$.
If otherwise, then $U$ is a proper subgroup of $W$.
By Maschke's Theorem, $V$ can be decomposed as
$W=U\times V$ such that $H$ normalises both $U$ and $V$.
By the definition of $U$, $H$ is fixed-point-free on $V$.
Let $k=\dim(V)$. Then $k<d$.
By the above paragraph, we conclude that
$2md$ divides $p^k-1$, contrary to our assumption.
This establishes the claim.
So $G$ is a primitive Frobenius group.

For $y$ defined in Construction~\ref{im}, let $z=y^{-1}y^{\ov\t^{\frac{p-1}{2}}}=x_1^2$.
Then $x_1$ belongs to $\l N,y\r$. So does $x_i$ for $1\leqslant i\leqslant d$.
It follows that $\l N,y\r=X$. Thus $\Ga$ is connected.
It furthermore implies that $\l\ov c\r$ belongs to $X$, and so
$\soc(X)=W\times\l\ov c\r$.

Let $\s_i=\t_i^{\frac{p-1}{2}}$ where $1\leqslant i\leqslant d$.
Finally, as $\s_i^y=\s_{i-1}$ for $3\leqslant i\leqslant d$, $\s_1^y=\s_d$ and $\s_2^y=x_1^2\s_1$,
we obtain that $N\cap N^y=\l\s_2,\s_3,\ldots,\s_d\r\cong\ZZ_2^{d-1}$.
That is to  say, $N\cap N^y$ has index $2$ in $N$.
Since $X\leqslant\Aut\Ga$, $\Ga$  is not a cycle.
By Lemma~\ref{Cos}, $\Ga$ is connected, $X$-edge-transitive and of valency 4.
\qed

\vskip0.07in
\noindent{\bf Remark.} The normal quotient $\Ga_W$ induced by $W$ is a cycle (see Lemmas~\ref{Lee} and \ref{Le}).

\vskip0.1in
As a matter of fact, there are several groups which are primitive Frobenius groups and satisfy Construction~\ref{im}.
For example, $G=\ZZ_7^3{:}\ZZ_{18}$, $\ZZ_{13}^3{:}\ZZ_{18}$, $\ZZ_{41}^5{:}\ZZ_{50}$, and so on.
\vskip0.1in

The following construction produces edge-transitive graphs admitting
a group $X$ satisfying part~(4) of Theorem~\ref{soluble} with $O=1$.

\begin{construction}\label{yin}
Using the notation in Construction~$\ref{im}$. Assume $\ell\geqslant 2$.
Let $N=\prod_{i\not=3}\l\t_i^{\frac{p-1}{2}}\r\cong\ZZ_2^{d-1}$,
and $h=\t_1^{\frac{p-1}{4m}}(12\cdots d)$. Let $X=W{:}\l N,h\r$, and let $G=W{:}\l h\r$. Set $y=(x_2h)^{-1}$, and
\[\Ga(p,n)=\Cos(X,N,N\{y,y^{-1}\}N).\]
\end{construction}

\begin{lemma}
Let $\Ga=\Ga(p,n)$ be a graph constructed in Construction~$\ref{yin}$. Then $\Ga$
is a connected tetravalent $X$-edge-transitive Cayley graph of Frobenius group $G$,
and $G$ is not  normal in $X$.
\end{lemma}
\proof Obviously, $G$ is not normal in $X$. Let $\s_i=\t_i^{\frac{p-1}{2}}$ where $1\leqslant i\leqslant d$.
By easy calculations, $\s_1^y=\s_d$, $\s_2^y=\s_1$, and $\s_i^y=\s_{i-1}$ for $4\leqslant i\leqslant d$.
It follows that $\s_3$ belongs to $\l N,y\r$, and $N\cap N^y=\l\s_1,\s_4,\ldots,\s_d\r\cong\ZZ_2^{d-2}$.
At the same time, we obtain $\s_3^y=x_2^2\s_2$, and hence $x_2$ belongs to $\l N,y\r$.
So does $x_i$ for each $i$. It implies that $\l N,y\r=X$. Consequently, $\Ga$ is connected.
Arguing similarly as Lemma~\ref{Liu}, we can obtain
that $G$ is a Frobenius group, and
$\Ga$ is $X$-edge-transitive Cayley graph of $G$ and of valency 4,
the statement follows.
\qed

\vskip0.07in
\noindent{\bf Remark.} Clearly, $\l h\r$ does not normalise $N$. In other words,
$X$ can't satisfy the properties in part~(a) of Lemma~\ref{Lee}.
However, $h$ normalises $\l N,h^{\frac{n}{2}}\r$, namely,
$X$ satisfies the properties in part~(ii) of Lemma~\ref{Le}.
Thus $\Ga_W\cong\C_{\frac{n}{2}[2]}$, where $N$, $W$, and $\Ga$ appear in Construction~\ref{yin},
(see Lemma~\ref{Lee} and Lemma~\ref{Le}).

\vskip0.07in
The following construction produces edge-transitive graphs admitting
a group $X$ satisfying part~(5) of Theorem~\ref{soluble} with $O=1$.

\vskip0.07in

Let $n=3p_1^{l_1}\cdots p_s^{l_s}$ be an odd number,
where $3$, $p_1$, $\ldots$, $p_s$ are pairwise distinct primes,
and $l_i\geqslant1$ for each $i$.
Let $G_1=W_1{:}H_1\cong\ZZ_2^d{:}\ZZ_n$
be a primitive Frobenius group.
Let $G_2$ be a subgroup of $ G_1$ such that $G_2=W_2{:}H_2\cong\A_4$.
Write $H_1=\l h_1\r$, and $H_2=\l h_2\r$ where $h_2=h_1^{\frac{n}{3}}$.

Let $H=\l h\r$ where $h=(h_1,h_2)$. Let $V=W_1\times W_2$,
and $W=\{(w,1)\div w\in W_1\}$. Set \[\mbox{$G=W{:}H$\, and\,
$X=V{:}H$}.\]
By the definition, it is easy to show that $G$ is a primitive Frobenius group.

\begin{construction}\label{D{2p}}
{\rm  Let $R=\l (w,w), (w,w)^h\r$, where $1\not=w\in W_2$. Set
\[\Ga(2,d,n)=\Cos(X, R, R\{h,h^{-1}\}R).\]
}
\end{construction}

\begin{lemma}
Let $\Ga=\Ga(2,d,n)$ be a graph constructed in Construction~$\ref{D{2p}}$. Then $\Ga$
is a connected tetravalent $X$-edge-transitive Cayley graph of $G$,
and $G$ is not  normal in $X$.
In particular, $\Ga_W$ is a cycle.
\end{lemma}
\proof We first prove that $R$ is core-free in $X$.
By the definition of $R$, we have $R\cong\ZZ_2^2$.
Assume $K\leqslant R$, and $1\not=K\unlhd X$.
Then $\Aut(K)$ is isomorphic  to a subgroup of $\S_3$.
So we conclude that $\C_H(K)\not=1$, which  contradicts the fact that $G_1$ is a Frobenius group.
Thus $R$ is core-free in $X$.
We observe that $R\cap G=1$, it follows that $|X|=|R||G|$,
and so $X=RG$. It implies that $G$ is regular on the vertex set $[X{:}R]$,
and hence $\Ga$ is a Cayley graph of $G$.

By the definition of $G_1$, we conclude that $w^{h_1^3}\not=w$ for any $1\not=w\in W_1$.
It implies that $(w,w)^{h^3}(w,w)\not=1$, namely, $(w^{h_1^3}w,1)\not=1$.
Since $H$ is irreducible on $W_1$, implying that $W_1$ belongs to $\l R,h\r$.
So does $V$. Thus $\Ga$ is connected.
Arguing similarly as above,  $G$ is not normal in $X$.
Clearly, $R\cap R^h=\l (w,w)^h\r\cong\ZZ_2$.
It follows that $R\cap R^h$ has index $2$ in $R$.
As $X\leqslant\Aut\Ga$, $\Ga$ is not a cycle, so
that by Lemma~\ref{Cos}, we obtain $\Ga$ is $X$-edge-transitive and of valency 4.
Note that $\Ga_W$ is a Cayley graph.
By \cite[Theorem~1.2]{Baik}, we conclude that $\Ga_W$ is a cycle.
\qed

By Constructions~\ref{BBB}-\ref{yin}, each case of Theorem~\ref{soluble} occurs.

\vskip0.07in
We now construct some examples of graphs appearing in Theorem~$\ref{isoluble}$.

\vskip0.07in
Based on several previous known results, arc-transitive elementary abelian covers of the complete
graph $\K_5$ were classified by Bo$\check{s}$tjan Kuzman \cite{Kuzman}. However,
for the completeness, we present here a distinct and independent construction.

\vskip0.07in
Let $p$ be a prime such that $5$ is a primitive divisor of $p^4-1$.
Set \[V=\l e_1\r\times\cdots\times\l e_5\r\cong\ZZ_p^5.\]
We define an action of $\A_5$ on $V$ as follows:
\[\mbox{$(\prod_{i=1}^5e_i^{\lambda_i})^{g}=\prod_{i=1}^5 e_{i^{g^{-1}}}^{\lambda_i}$,
where $g\in\A_5$, and $0\leqslant\lambda_i\leqslant 4$ for each $i$.}\]
By this definition, $\A_5$ acts naturally on $V$.
Let  $\ov e_i=e_5e_i^{-1}$ for $1\leqslant i\leqslant 4$.
Set
\[W=\l\ov e_1\r\times\l\ov e_2\r\times\l\ov e_3\r\times\l\ov e_4\r.\]
It is straightforward to show that $\A_5$ acts faithfully on $W$.

\begin{construction}\label{AA}
{\rm  Let $G=W{:}\l h\r$ with $h=(12345)$, and
let $X=W{:}N=\ZZ_p^4{:}\A_5$.
Let $R=\Alt\{2,3,4,5\}\cong\A_4$, and let $g=\ov e_1(15)(24)$. Set
\[\Ga(p,4,5)=\Cos(X, R, RgR).\]
}
\end{construction}

\begin{lemma}\label{A_5}
Let $\Ga=\Ga(p,4,5)$ be a graph constructed in Construction~$\ref{AA}$.
Then $\Ga$ is a connected tetravalent $(X,2)$-arc-transitive Cayley graph of Frobenius group $G$,
and $G$ is not normal in $X$.
In particular, $\Ga$ is a cover of $\Ga_W$, and $\Ga_W\cong\K_5$.
\end{lemma}
\proof Let $H=\l h\r$.
By definition of $W$, we conclude $H$ is fixed-point-free on $W$.
Since $5$ is a primitive divisor of $p^4-1$,
$H$ acts irreducibly on $W$.
That is to say, $G$ is a primitive Frobenius group.
Clearly, $N$ has a decomposition $HR$.
It implies that $R$ is core-free in $X$, and hence $X\leqslant\Aut\Ga$.
Now $X=GR$ and $G\cap R=1$, and so $G$ is regular on the vertex set $[X{:}R]$.
Thus $\Ga$ is a Cayley graph of $G$.
Obviously, $G$ is not normal in $X$.

Denote by $u$ and $v$ the vertices $R$ and $Rg$, respectively. Then $X_u=R$ and $X_v=R^g$.
Let $r=(234)$.
Since $X_{uv}=X_u\cap X_v$, a small calculations show $X_{uv}=\l r\r$.
By Lemma~\ref{Cos}, $\Ga$ is of valency $4$.
It is clear that $g$ has order $2$, and $r^g=r^{-1}$.
So $g\in\N_X(X_{uv})$.
Let $\ov R=\l R,g\r$.
Since $(15)(24)=(12345)(25)(34)$, we conclude that $\ov e_1h$ belongs to $\ov R$.
Let $a=(25)(34)$ and $b=(23)(45)$.
By easy calculations, we obtain
\[\mbox{$\ov e_1h(\ov e_1h)^a=\ov e_1\ov e_2\ov e_3^{-1}$,
$(\ov e_1\ov e_2\ov e_3^{-1})^b=\ov e_1\ov e_2^{-1}\ov e_3\ov e_4^{-1}$, and
$(\ov e_1\ov e_2^{-1}\ov e_3\ov e_4^{-1})^{ab}=\ov e_1\ov e_2^{-1}\ov e_4^{-1}$}.\]
Combining the above three equations, we conclude that $\ov e_3$ belongs to $\ov R$.
So does $\ov e_i$ for $i=1,2$, $4$.
Consequently, $W\leqslant\ov R$.
Recall that $\ov e_1h$ is inside in $\ov R$,
it follows that $h$ belongs to $\ov R$,
forcing $\l R,g\r=X$.
Thus $\Ga$ is connected.
Since $X/W$ is insoluble,  by Lemma~\ref{insoluble},
$\Ga$ is a cover of $\Ga_W$. Clearly, $\Ga_W$ is a Cayley graph of $G/W$.
By \cite[Theorem~1.2]{Baik}, we obtain $\Ga_W\cong\K_5$.
\qed

Let $p$ be a prime such that $10$ is a primitive divisor of $p^4-1$.
Let
\[
V_1=\l e_1\r\times\cdots\times\l e_5\r,\, \mbox{and}\,\,\,
V_2=\l e_{1'}\r\times\cdots\times\l e_{5'}\r
\]
such that $V_1\cong V_2\cong\ZZ_p^5$.
Set $T=\l(12345)(1'2'3'4'5'),(12)(1'2')\r$.
It is straightforward to show that $T\cong\S_5$.
Then, for any $g\in T$,  $g$ acts on $V_i(i=1,2)$ as follows:
\[\begin{array}{l}
(\prod_{i=1}^5e_i^{\lambda_i})^{g}=\prod_{i=1}^5{e_{{i}^{g^{-1}}}^{\lambda_i}}, \mbox{where}\,\ 0\leqslant \lambda_i\leqslant 4\,\, \mbox{for each}\,\, i,\\
 (\prod_{i=1}^5e_{i'}^{\lambda_{i'}})^{g}=\prod_{i=1}^5e_{{i'}^{g^{-1}}}^{\lambda_{i'}}, \mbox{where}\,\ 0\leqslant \lambda_{i'}\leqslant 4\,\, \mbox{for each}\,\, i'.
 \end{array}\]

Let $\ov e=\prod_{i=1}^5e_i$, and $\ov e'=\prod_{i=1}^5e_{i'}$.
Let  $\ov e_i= e_i\l\ov e\r$ and $\ov e_{i'}=e_{i'}\l \ov e'\r$ for $1\leqslant i\leqslant 5$.
Set
\[\mbox{$W=\l w_1\r\times\l w_2\r\times\l w_3\r\times\l w_4\r$ where $w_i=\ov e_i\ov e_{i'}^{-1}$.}\]
Then $W\cong\ZZ_p^4$.
Note that $T$ fixes each element of $\l\ov e\r$ and $\l\ov e'\r$.
So $T$ induces a faithful action on $W$.
Without loss of generality, we may assume that $T$ is a subgroup of $\GL(W)$.
Let $g=(11')\cdots(55')$. Obviously, $g$ inverts each non-identity element of $W$.

\begin{construction}\label{AAAAA}
{\rm Let $G=W{:}H$ where $H=\l h, g\r$ with $h=(12345)(1'2'3'4'5')$.
Set $X=W{:}N=W{:}(T\times \l g\r)\cong\ZZ_p^4{:}(\S_5\times\ZZ_2)$. %where $5$ is a primitive of $p^4-1$.
%Assume that $G=W{:}H$ is a primitive Frobenius group satisfying that
Let $R=\l (1234)(1'2'3'4'),(12)(1'2')\r\cong\S_4$,
and let $y=w_1w_5(15)(1'5')g$. Set
\[\Ga(p,4,10)=\Cos(X, R,RyR).\]
}
\end{construction}

Arguing similarly as Lemma~\ref{A_5}, we have the following conclusion in next lemma.
\vskip0.07in
\begin{lemma}

Let $\Ga=\Ga(p,4,10)$ be a graph constructed in Construction~$\ref{AAAAA}$.
Then $\Ga$ is a connected tetravalent $(X,2)$-arc-transitive Cayley graph of Frobenius group $G$,
and $G$ is not normal in $X$.
In particular, $\Ga$ is a cover of $\Ga_{W}$, and $\Ga_W\cong\K_{5,5}-5\K_2$.
\end{lemma}

\vskip0.07in

Here are a few of graphs whose automorphism groups are almost simple.

\begin{example}\label{TA}
{\rm Let $\FF=\GF(p)$ be a finite field of order $p$.
Let $U$ and $V$ consist of 1-subspaces and 2-subspaces of $\FF^3$, respectively.

\vskip0.07in
{\bf Case 1{:}}\
Let $p=2$. Define a bipartite graph $\Ga$ with bipartite $U$ and $V$ such that $u\in U$ and $v\in V$ are adjacent if and only if
$u+v=\FF^3$. This is the point-line non-incidence graph of the Fano plane $\PG(2,2)$.
Furthermore, $\Aut\Ga=\PGL(3, 2).\ZZ_2$, and $\Ga$ is
a Cayley graph of $G=\D_{14}$. For example, refer to \cite{p}.

\vskip0.07in
{\bf Case 2{:}}\
Let $p=3$. Define a bipartite graph $\Ga$ with bipartite $U$ and $V$ such that $u\in U$ and $v\in V$ are adjacent if and only if $u$ is
a subspace of $v$. Then $\Ga$ is the point-line incidence graph of the projective plane $\PG(2,3)$.
Furthermore, $\Aut\Ga=\PGL(3,3).\ZZ_2$,
and $\Ga$ is a Cayley graph of $G=\D_{26}$. Refer to \cite{Li,flag}, for example.

}
\end{example}

\begin{example}\label{A}
{\rm
Let $X=\PGL(2,7)$. By the Atlas \cite{atlas}, $X$ has a maximal subgroup $H\cong\D_{16}$. Pick a subgroup $K\leqslant H$
with $K\cong\ZZ_2^2$. Then $\D_8\cong\N_{H}(K)\leqslant\N_{X}(K)\cong\S_4$. Choose an involution $o\in\N_{H}(K)\setminus K$
and an element $z\in\N_X(K)$ of order $3$ such that $z^o=z^{-1}$. Then $\langle o, z\rangle\cong\S_3$,
and $o(oz)=2$. Since $H$ is a maximal subgroup of $X$, it follows that $\langle H, oz\rangle=X$.
Let $\Ga=\Cos(X, H, HozH)$. By the choices of $o$ and $z$, we conclude that $|H{:}H\cap H^{oz}|=4$, namely,
$\Ga$ is a connected $X$-arc-transitive graph of valency $4$.
By MAGMA~\cite{Magma},  we have that $X=GH$ where $G=\ZZ_7{:}\ZZ_3$,
and thus $\Ga$ is a connected $X$-arc-transitive Cayley graph of $G$ of valency $4$.
By Li et al.\cite{Zai Ping Lu}, $\Aut\Ga=\PGL(2,7)$.
}
\end{example}

\begin{example}\label{B}
{\rm
Let $X=\PGL(2,7)$. Then $T=\soc(X)\cong\PSL(2,7)$.
Take $H\leqslant T$ such that $H\cong\D_8$.
Choose an involution $o$ such that $o$ is not in the center of $H$.
It is simple to check that
$\N_{H}(\langle o\rangle)\cong\ZZ_2^2$, $\N_{T}(\langle o\rangle)\cong\D_8$
and $\N_X(\langle o\rangle)\cong\D_{16}$.
Let $\N_{X}(\langle o\rangle)=\N_{T}(\langle o\rangle){:}\langle z\rangle$
for some involution $z\in X\setminus T$.
Take $y\in\N_{T}(\langle o\rangle){:}\langle z\rangle$  of order $4$.
Set $\Ga=\Cos(X, H, HxH)$, where $x=z$ or $yz$.
By Li et al.\cite{Zai Ping Lu}, $\Ga$ is a connected tetravalent arc-transitive Cayley graph of $\ZZ_7{:}\ZZ_6$,
and $\Aut\Ga=\PGL(2,7)$.
}
\end{example}

\begin{example}\label{C}
{\rm
Let $X=\PGL(2,11)$ or $\PSL(2,23)$. By
the Atlas \cite{atlas}, $X$ has a maximal subgroup $H\cong\S_4$.
Let $L\cong\S_3$ be a subgroup of $H$.
Checking the subgroups of $X$ in the Atlas \cite{atlas},
we conclude that $\N_X(L)=\l o\r\times L\cong\D_{12}$,
where $o\in \N_X(L)\setminus H$ is an involution.  Set $\Ga=\Cos(X,H,HoH)$. Since $H$
is a maximal subgroup of $X$, $\langle o, H\rangle=X$. It is straightforward to check that $|H{:}H\cap H^o|=4$.
Then $\Ga$ is a connected tetravalent $X$-arc-transitive graph.
Moreover, $X$ has a subgroup $G$ which is regular on the vertices, where $G\cong\ZZ_{11}{:}\ZZ_5$
or $\ZZ_{23}{:}\ZZ_{11}$, respectively. We denote by $P_{11,5}$ and $P_{23,11}$ the graphs
associated with $\PGL(2,11)$ and $\PSL(2,23)$, respectively.
By Li et al.\cite{Zai Ping Lu},
$\Aut P_{11,5} =\PGL(2,11)$ and $\Aut P_{23,11} =\PSL(2,23)$.
}
\end{example}

\begin{example}\label{D}
{\rm
Let $\Ga=(V\Ga, E\Ga)$ be a connected arc-transitive Cayley graph. The standard double cover $\Ga^{(2)}$
is the graph with vertex set $V\Ga\cup\{u'|u\in V\Ga\}$ such that $\{u,v'\}\in E\Ga^{(2)}$ whenever
$\{u,v\}\in E\Ga$.  For each $x\in\Aut\Ga$, define $\tilde x : u\to u^x$, $u'\to (u^x)'$.
Then $\Aut\Ga$ can be viewed as a subgroup of $\Aut\Ga^{(2)}$ in this way.
Define $\epsilon : u\to u'$, $u'\to u$. Then $\epsilon\in\Aut\Ga^{(2)}$.
Set $X=\l\Aut\Ga,\epsilon\r$. Then $X=\Aut\Ga\times\l\epsilon\r$. So
$\Ga^{(2)}$ is an $X$-arc-transitive Cayley graph. By Li et al.\cite{Zai Ping Lu},
$P_{11,5}^{(2)}$ is a Cayley graph of $\ZZ_{11}{:}\ZZ_{10}$, and $\Aut P^{(2)}_{11,5}\cong\PGL(2,11)\times\ZZ_2$.
}
\end{example}

\section{Soluble automorphism groups}
In this section, let $G=W{:}H\cong\ZZ_p^d{:}\ZZ_n$ be a primitive Frobenius group.
Let $\Ga=\Cay(G,S)$ be a connected  $X$-edge-transitive tetravalent Cayley graph,
where $G\leqslant X\leqslant \Aut\Ga$.
Denote by $F$  the Fitting subgroup of $X$.
If $X$ is solvable,
then an important property of its Fitting subgroup is that it is self-centralized,
that is, $\C_X(F)\leqslant F$.
In what follows, we will determine the graph $\Ga$ for the case where $X$ is solvable.

\begin{lemma}\label{un}
Assume that $F$ is a $r$-group, where $r$ is a prime.
If $\Ga$ is a cover of $\Ga_F$,
then $F=W$.
\end{lemma}
\proof
Note that $W$ is minimal and normal in $G$.
Then either $W\leqslant F$ or $F\cap G=1$.
If $W\leqslant F$, then $F=W$ as $\Ga$ is a cover of $\Ga_F$.
Thus we assume that $F\cap G=1$.
%Suppose, to the contrary, that $F\cap G=1$.

Let $\ov G=G\Phi(F)/\Phi(F)$, and let $\ov F=F/\Phi(F)$.
Since $\Phi(F)$ \char $F$, we obtain $\ov G$ can act on $\ov F$ by conjugation.
Clearly, $\ov G\cong G$. In what follows, write $\ov G=\ov W{:}\ov H\cong W{:}H$.

Assume first that $r\not=p$.  If $\ov W$ acts trivially on $\ov F$, then
$W$ induces the identity on $\ov F$.
From \cite[p.174, Theorem~1.4]{Gorenstein},
it follows that $W$ acts trivially on $F$.
So $W\leqslant\C_X(F)\leqslant F$, against our assumption.
Thus $\ov W$ acts nontrivially on $\ov F$.

Let $M=\ov F{:}\ov G$.
Let \[1\unlhd M_1\unlhd M_2\unlhd\cdots\unlhd M_{m-1}\unlhd M_m=\ov F\]
be the normal series of $\overline{F}$ such that each
$M_i/M_{i-1}$ is a minimal normal subgroup of $M/M_{i-1}$ for $1\leqslant i\leqslant m$,
where $M_0=1$.

It is straightforward to show that $\ov G$ normalizes $\C_{M_1}(\ov W)$,
and hence $\C_{M_1}(\ov W)\unlhd M$.
By the minimality of $M_1$,
we conclude that either $\C_{M_1}(\ov W)=1$ or $\C_{M_1}(\ov W)=M_1$.
If  $\C_{M_1}(\ov W)=1$, so that
by \cite[Theorem~2.7]{Evgenii I} we have $|M_1|=|\C_{M_1}(\ov H)|^{|\ov H|}$.
However,  $\Ga$ is a normal cover of $\Ga_F$, we conclude that $|M_1|$ divides $|\ov H|$, a contradiction occurs.
Thus $\C_{M_1}(\ov W)=M_1$, that is, $\ov W\leqslant\C_M(M_1)$.
%Arguing as above with $M_2/M_1$ in the place $M_1$, we obtain that
%$\ov WM_1/M_1\leqslant\C_{M/M_1}(M_2/M_1)$,
%namely, $[M_2,\ov W]\subseteq M_1$.
Repeating the above argument for $M_{i}/M_{i-1}$ and $\ov G M_{i-1}/M_{i-1}$,
we obtain that $[M_i,\ov W]\subseteq M_{i-1}$ for $1\leqslant i\leqslant m$.
It follows that $\ov W$ stabilizes the normal series of $\ov F$,
and hence $\ov W$ centralizes $\ov F$, refer to \cite[p.178, Theorem~3.2]{Gorenstein}.
It implies that $W$ induces the identity on $\ov F$,
and so $W$ is trivial on $F$ (see \cite[p.174, Theorem~1.4]{Gorenstein}), namely,  $W\leqslant\C_X(F)$, again
against our assumption.

Assume now that $r=p$.
Then $|F|\leqslant|W|$.
Denote by $\Sigma$ the normal quotient graph $\Ga_F$.
If $|F|=|W|$,  then $W$ fixes each vertex of $\Sigma$,
and hence $W\leqslant F$, which is impossible.
Thus $|F|<|W|$.
Set $\ov X=X/F$.
Let $F_{\ov X}$ be the Fitting subgroup of $\ov X$.
It is known that $F_{\ov X}$ is a $p'$-group.
Let $\tilde G=GF/F\cong G$. Write $\tilde G=\tilde W{:}\tilde H$.
For that case, we conclude $F_{\ov X}\cap \tilde G=1$.
It follows that $\Ga$ is $(X,2)$-arc-transitive, and so
$\Sigma$ is $(\ov X,2)$-arc-transitive.
%Now we consider  the normal quotient graph $\Sigma_{F_{\ov X}}$.
By \cite[Theorem~4.1]{CE}, $\Sigma$ is a cover of $\Sigma_{F_{\ov X}}$ or $\Sigma_{F_{\ov X}}=K_2$.
For the former, arguing as above,
we also obtain $\tilde W\leqslant F_{\ov X}$,
which contradicts $F_{\ov X}\cap \tilde G=1$.
For the latter, we obtain $p=2$, and $|G|$ divides $2^53^6$.
Since $G$ is a primitive Frobenius group, we have $G\cong\ZZ_2^2{:}\ZZ_3$.
So $F\cong\ZZ_2$, and thus $F\leqslant\Z(X)$, again a contradiction.
Therefore, $F=W$.
\qed

For a group $T$ and a prime $q$,
by $T_q$ we mean a Sylow $q$-subgroup of $T$.

\begin{lemma}\label{Lei}
Use the notation defined above, rewrite $\Ga=(V\Ga,E\Ga)$.
Then we have:
\begin{itemize}
\item[(i)]
If $p$ is an odd prime, then either $G\cong\D_{2p}$ or $W\unlhd X$;

\item[(ii)]
If $p=2$, then $F=O_2(X)$, and
\begin{itemize}
\item[(a)]
$W<F$, $\Ga_F$ is a cycle, and $X=(F{:}H).\calO$,
where $\calO=1$ or $\ZZ_2$;
\item[(b)]
$W=F$ and $W\unlhd X$.
\end{itemize}
\end{itemize}
\end{lemma}
\proof Suppose that $G\ncong\D_{2p}$. We first claim that $W\leqslant F$ and $F\cap H=1$.
Suppose, by way of contradiction, that $G\cap F=1$.
Since $X=GX_1$, we obtain  $|F|$ divides $|X_1|$. From  Lemma~$\ref{Lu}$, it follows that
each prime divisor of $|F|$ is either $2$ or $3$.

Let $K$ be the kernel of $X$ acting on $\Ga_F$.
Then $X/K\leqslant\Aut\Ga_F$.
Recall that $p$ is a prime divisor of $|W|$.
Suppose that $p>3$.
Let $B$ be an orbit of $F$ acting on $V\Ga$.
So $|B|$ divides $|F|$.
If $G\cap K\not=1$, then $W\leqslant K$.
Let $\bigtriangleup$ be an orbit of $W$ acting
on $V\Ga$, which is contained in the block $B$.
Then $|\bigtriangleup|$ divides $|B|$,
which is impossible.
So $G\cap K=1$. Let $\ov G=GK/K$.
Then $\ov G\cong G$ is a Frobenius group. Write $\ov G=\ov W{:}\ov H\cong G$.
If $\Ga_K$ is a cycle, then $d=1$, and so $\ov G\cong\D_{2p}$, against our assumption.
Thus $\Ga$ is a cover of $\Ga_K$, and then $K=F$.
By Lemma~\ref{Lu}, $\ov G_B\leqslant \ov H$.
In view of Lemma~$\ref{un}$, $F$ is a $\{2,3\}$-group.
So is $\ov G_B$ because $|\ov G_B|=|F|$.
Note that $\Ga_F$ is $\ov G$-vertex-transitive.
By \cite[Lemma~2.1]{Cai heng},  we conclude that $\ov G_B^{\Ga(B)}$
is a cyclic  group of order $2^i3^j$ for $i,j\geqslant1$.
However, $\ov G_B^{\Ga(B)}$ is isomorphic to a subgroup of $\S_4$,
which is a contradiction. Thus $p=2$ or $3$.

If $F=O_2(X)$,
by Lemma~$\ref{un}$, $\Ga_F$ is a cycle.
So is $\Ga_K$.
By the assumption, we conclude that $p=2$, and $W\leqslant K$.
It follows that $K=W.K_1$, where $K_1$ is a $2$-group.
Since  $K\unlhd X$, we have $K\leqslant F$, which contradicts the fact that $F\cap G=1$.

If $F=O_3(X)$, then $\Ga$ is $(X,2)$-arc-transitive.
From Lemma~\ref{un}, it follows that $\Ga_F=K_2$.
It implies that  $G\cong\D_6$,
against our assumption.

If  $F=O_2(X)\times O_3(X)$, then  $\Ga$ is $(X,2)$-arc-transitive.
By \cite[Theorem~4.1]{CE}, $\Ga$ is a cover of $\Ga_F$, $\Ga_F=K_2$ or $F$ is transitive on $V\Ga$.
Assume first that $\Ga$ is a cover of $\Ga_F$.
Let $Y=\Aut\Ga_F$ and $\ov X=X/F$.
Then $\ov X$ is a subgroup of $Y$.
Recall that $F\cap G=1$, we have $|\ov G_{B}|=|F|$, and hence
$\ov G_{B}$ is a Frobenius group of $Y_{B}$.
Since $\Ga_F$ is $\ov G$-vertex-transitive,
we conclude that $\ov G_{B}^{\Ga_F(B)}=\A_4$, $\S_3$ or $\S_4$, refer to \cite[Lemma~2.1]{Cai heng}.
Assume $\ov G_{B}^{\Ga_F(B)}=\A_4$ or $\S_4$.
Then $\Ga_F$ is $(\ov G,2)$-arc-transitive.
From Lemma~\ref{Lu}, it follows that
$\ov G_{B}=\A_4$, which implies that $|F|=12$.
For this case, it is easy to show that $W\leqslant\C_X(F)\leqslant F$, a contradiction occurs.
Assume $\ov G_B^{\Ga_F(B)}=\S_3$. Let $\ov G_B^{[1]}$ be the kernel of $\ov G_B$ acting on $\Ga_F(B)$.
Recall that $\ov G_{B}$ is a Frobenius group,
we conclude $\ov G_B^{[1]}$ is a $3$-group, and hence $O_2(X)\cong\ZZ_2$.
It further implies that $W$ is a $3$-group.
Let $\ov F$ be the Fitting subgroup of $\ov X$.
Then $\ov F\cap \ov G=1$. It follows that $\ov F$ is a $2$-group.
Since $|F.\ov F|$ divides $2^43^6$, we have $|\ov F|\leqslant 8$,
and thus  $\ov W\leqslant\C_{\ov X}(\ov F)$, again a contradiction.

Assume now that $\Ga_F=K_2$. Then  $|V\Ga|$ divides $2^53^6$.
Since $\Ga$ is a cover of $\Ga_{F_2}$ and $\Ga_{F_3}$,
we conclude that $|F_2|=\frac{|G_2|}{2}$ or $|G_2|$, and $|F_3|=|G_3|$.
When $p=3$,  we have  $|F_3|=|W|$.
Note that $W$ is minimal and normal in $G$.
Then $W$ fixes each vertex of $\Ga_{F_3}$, and thus $W\leqslant F$, which contradicts $G\cap F=1$.
For $p=2$, and $|F_2|=|G_2|$,  we also obtain the same contradiction.
When $p=2$ and $|F_2|=\frac{|G_2|}{2}$.
Since $G$ is a $\{2,3\}$-group, we conclude that $G\cong\ZZ_2^2{:}\ZZ_3$,  and hence $|F|=6$.
For that case, we easily obtain that $W\leqslant F$, again a contradiction.
Similarly, we also exclude the case where $F$ is transitive on $V\Ga$.

Summarizing the above discussion, we obtain $W\leqslant F$.
Since $G$ is a Frobenius group, we have $F\cap H=1$, as claimed. Next
we process our analysis by several cases.

\vskip0.07in
{\bf Case~1:}\ If $p>3$, then  $W\unlhd X$.

By the previous discussion, we have $W\leqslant F_p$. By Lemma~\ref{Lu}, we conclude that $W=F_p$,
and hence $W\unlhd X$.

\vskip0.07in
{\bf Case~2:}\ If $p=3$, then $W\unlhd X$.

 If $W<F_3$, then $\Ga$
is $(X,2)$-arc transitive.
For this case, $\Ga_{F_3}=K_2$, and so $G\cong\D_6$,
contrary to our assumption. Thus $W=F_3$, and then $W\unlhd X$.

\vskip0.07in
{\bf Case~3:}\ If  $p=2$, then either $\Ga_F$ is a cycle, or $W\unlhd X$.

Assume that $W<F_2$. Since $F\cap H=1$, we know that $|FH|=|F||H|$.
Note that $|FH|$ divides $|X|$, it follows that $|F|$ divides $|W||X_1|$,
and hence $|F_{2'}|$ divides $|X_1|$.
So $F_{2'}$ is a $3$-group. If $F_{2'}\not=1$, then $\Ga$ is $(X,2)$-transitive, and thus it follows from \cite[Theorem~4.1]{CE}
that $\Ga$ is a cover of $\Ga_{F_2}$, a contradiction.
So $F_{2'}=1$. By Lemma~\ref{un}, $\Ga_F$ is a cycle because $|H|$ is an odd number.
Recall that $B$ is a vertex of $\Ga_F$.
Since $\Ga$ is a Cayley graph of $G$,
we obtain $W$ is regular on $B$.
So $K=F=WK_1$.
Consequently, $X=(F{:}H).\calO$ where $\calO\cong1$ or $\ZZ_2$.
This completes the proof of Lemma~\ref{Lei}.
\qed

For the group $G\cong\D_{2p}$ where $p$ is an odd prime.
Applying \cite[Theorem~1.1]{Zai Ping Lu},
we have the following conclusions.

\begin{lemma}\label{D_2p}
Let $G\cong\D_{2p}$, and let $\Ga$ be a connected edge-transitive tetravalent Cayley graph of $G$.
Then we have
\begin{itemize}
\item[(i)]
$\Ga$ is arc-regular, and $\Aut\Ga\cong\D_{2p}{:}\ZZ_4$;
\item[(ii)]
$\Ga\cong\C_{p[2]}$, and $\Aut\Ga\cong\ZZ_2^p{:}\D_{2p}$.
\end{itemize}
\end{lemma}

In the remainder of this section assume
that $G\not\cong \D_{2p}$ with $p$ an odd prime, unless specified otherwise.

\vskip0.01in

Recall that the {\it socle} of a finite group $R$ (denoted by $\soc(R)$) is the product of
all minimal normal subgroups of $R$. Evidently, $\soc(R)$ is a characteristic subgroup of $R$.

\vskip0.01in

We next treat the case where $W\unlhd X$, and
the normal quotient $\Ga_W$ is a cycle.

\begin{lemma}\label{Lee}
Let $K$ be the kernel of $X$ acting on $\Ga_W$.
Then the following statements hold:
\begin{itemize}
\item[(i)]
$X=((WK_1){:}H).\calO$, and $W\cong\ZZ_2^d$,
where $\calO\cong1$ or $\ZZ_2$;

\item[(ii)]
Assume $p$ is an odd prime. Then we have
\begin{itemize}
\item[(1)]
 $G$ is normal in $ X$, or
 \item[(2)]
 $G$ is not normal in $X$, and
\begin{itemize}
\item[(a)]
$X=W{:}((K_1{:}H).\calO)$, and $H$ does not centralise $K_1$
where $K_1\cong\ZZ_2^l$ with $2\leqslant l\leqslant d$, and $\calO\cong1$ or $\ZZ_2$;
\item[(b)]
 there exist $x_1,x_2,\ldots,x_d\in W$ and $\t_1,\t_2,\ldots,\t_d\in K_1$ such that $W=\l x_1,\ldots,x_d\r$,
 $\l x_i,\t_i\r\cong\D_{2p}$, and $K_1=\l\t_i\r\times\C_{K_1}(x_i)$ for $1\leqslant i\leqslant d$;
 \item[(c)]
 $\soc(X)=W\times L$, where $L\cong 1$ or $\ZZ_2$;
 \item[(d)]
 $H$ is imprimitive on  $W$.
  \end{itemize}
\end{itemize}
\end{itemize}
\end{lemma}
\proof Let $B$ be a vertex of $\Ga_W$.
Since $\Ga$ is a Cayley graph of $G$,
we obtain $W$ is regular on $B$.
Thus $K=WK_1$, and $K\cap H=1$,
where $K_1$ is a $2$-group.
For that case,
$\Ga_W$ is a connected Cayley graph.
Recall that $H$ is of order $n$,
$\Ga_W$ is a cycle of size $n$, say.
It follows that $X/K\cong\ZZ_n$ or $\D_{2n}$. Further,
$\Ga$ is $X$-arc-transitive if and only if $X/K\cong\D_{2n}$.

Assume first that $p=2$.  Since $G\leqslant X$ and $(|K|,|H|)=1$,
we conclude that $K{:}H\leqslant X$.
Noting that $X/K$ is isomorphic to a subgroup of  $\D_{2n}$,
it follows that $X=(K{:}H).\calO$ with $\calO\cong1$ or $\ZZ_2$,
so we have part (i).

Assume now that $p$ is an odd prime.
Furthermore, we assume that $G$ is not normal in $X$.
If $K_1=1$, then $K=W$, and hence $G\lhd X$, which contradicts the assumption.
Thus $K_1\not=1$.

Let $U=\N_X(K_1)$. Since $K_1\ntrianglelefteq X$, it follows that $U\not=X$. Noticing that $(|W|, |K_1|) = 1$,
we obtain that $\N_{X/W}(K/W) =\N_{X/W}(WK_1/W) =\N_X(K_1)W/W = UW/W$. As $K/W\unlhd X/W$,
implying that $X = WU$. Since $W\lhd X$, we have that $W\cap U\lhd U$. Furthermore, $W\cap U\lhd W$ since $W$ is
abelian. Then $W\cap U\lhd\l U,W\r= UW=X$. If $W\leqslant U$, then $K=WK_1=W\times K_1$, and hence
$K_1\lhd X$, which is impossible. Thus $W\cap U<W$. Furthermore, note that $W$ is a minimal normal subgroup of $X$,
we obtain that $W\cap U=1$, and so $K\cap U=WK_1\cap U=(W\cap U)K_1=K_1$.
Now $X/K=UW/K=UK/K\cong U/(K\cap U)=U/K_1$, and hence $U=(K_1.\hat H).\calO$, where
$\hat H\cong\ZZ_n$ and $\calO\cong1$ or $\ZZ_2$.
Noting that $G$ belongs to $X$ and $G=W{:}H$,
there exists some $H^z\leqslant U$ such that $G=W{:}H^z$ is regular on $V\Ga$, where $z\in W$.
Without loss of generality, we may assume that  $U=(K_1{:}H).\calO$.
Then $X_1=K_1.\calO$. Furthermore, since $G$ is not
normal in $X$, we conclude that $H$ does not centralise $K_1$.

Set $Y=W{:}(K_1{:}H)$. Then $Y$ has index at most $2$ in $X$, and $\Ga$ is $Y$-edge-transitive.
It is obvious that $\Ga$ is not $Y$-arc-transitive.
Hence $\Ga=\Cos(Y,K_1,K_1\{y, y^{-1}\}K_1)$, where
$y\in Y$ is such that $\l K_1,y\r=Y$ and $K_1\cap K_1^y$
has index $2$ in $K_1$. We may choose $y\in W{:}H=G$
such that $H=\l h\r$ and $y=hx$ where $x\in W$. Then $K_1\cap K_1^y=K_1\cap K_1^x$
has index $2$ in $K_1$.

We claim that $K_1\cap K_1^x=\C_{K_1}(x)$.
For any $\sigma\in K_1\cap K_1^x$, we have that $\sigma^{x^{-1}}\in K_1$, and hence $\sigma^{-1}\sigma^{x^{-1}}\in K_1$.
Since $x\in W$ and $W\lhd WK_1$, we obtain that $\sigma^{-1}\sigma^{x^{-1}}=(\sigma^{-1}x\sigma)x^{-1}\in W$.
So $\sigma^{-1}\sigma^{x^{-1}}\in W\cap K_1=1$, and then $\sigma^{x^{-1}}=\sigma$.
Thus $\sigma$ centralises $x$.
It follows that $K_1\cap K_1^x\leqslant \C_{K_1}(x)$.
Clearly, $\C_{K_1}(x)\leqslant K_1\cap K_1^x$.
So $\C_{K_1}(x)=K_1\cap K_1^x$ as required.

Recall that $W$ is a minimal normal subgroup of $X$ and $X=WU$, we obtain that
$W=\l x\r\times\l x^{\sigma_2}\r\times\cdots\times\l x^{\sigma_d}\r$
where $\sigma_i\in U$. Then $\C_{K_1}(x^{\sigma_i})=(\C_{K_1}(x))^{\sigma_i}<K_1^{\sigma_i}=K_1$.
The intersection $\cap_{i=1}^d\C_{K_1}(x^{\sigma_i})\leqslant\C_K(W)=W$, and hence $\cap_{i=1}^d\C_{K_1}(x^{\sigma_i})=1$.
Since each $\C_{K_1}(x^{\sigma_i})$ is a
maximal subgroup of $K_1$, the Frattini subgroup $\Phi(K_1)\leqslant\cap_{i=1}^d\C_{K_1}(x^{\sigma_i})=1$.
Hence $K_1$ is an elementary abelian $2$-group, that is, $K_1\cong\ZZ_2^l$ for some $l\geqslant1$.
Recall that $\cap_{i=1}^d\C_{K_1}(x^{\sigma_i})=1$,
it follows that $l\leqslant d$.
Assume that $l=1$. Then $K_1\cong\ZZ_2$ and so
$K_1\leqslant\C_X(H)$. Thus $G\lhd X$,
which contradicts the fact that $G$ is not normal in $X$.
Hence $l>1$,   as in part (a).

Since $\C_{K_1}(x)$ has index $2$ in $K_1$, there exists some $\tau_1$
belonging to $ K_1$ such that $K_1=\l \tau_1\r\times\C_{K_1}(x)$.
Set $x_1=x^{-1}x^{\tau_1}$. Then $x_1\not=1$, $x_1^{\tau_1}=x_1^{-1}$ and $\C_{K_1}(x)=\C_{K_1}(x_1)$,
and hence $K_1=\l \tau_1\r\times\C_{K_1}(x_1)$.
Noticing that $W$ is a minimal normal subgroup of $X=WU$, there exist $\mu_1=1,\mu_2,\ldots,\mu_d\in U$
such that $W=\l x^{\mu_1}\r\times \l x^{\mu_2}\r\times\cdots\times\l x^{\mu_d}\r$.
Let $x_i=x_1^{\mu_i}$ and $\tau_i=\tau_1^{\mu_i}$, where $i=1,\ldots,d$.
Then $\ZZ_2^{l-1}\cong(\C_{K_1}(x_1))^{\mu_i}=\C_{K_1^{\mu_i}}(x_1^{\mu_i})=\C_{K_1}(x_i)$,
and $K_1=K_1^{\mu_i}=\l\tau_i\r\times\C_{K_1}(x_i)$. Furthermore,
$x_i^{\tau_i}=x_1^{\tau_1\mu_i}=(x_1^{-1})^{\mu_i}=x_i^{-1}$, and so $\l x_i,\tau_i\r\cong\D_{2p}$,
as in part (b).

Recall that $W\cong\ZZ_p^d$ for an odd prime $p$. Since $G$ is not normal in $X$, we conclude that $d>1$.
Assume that $X$ has a minimal normal subgroup $L\not=W$. Then $W\cap L=1$,
and $LK/K\lhd X/K\leqslant\D_{2n}$. It follows that
either $L\leqslant K$, or $L\cap K=1$. If $L\leqslant K$, then $L$ is a $2$-group.
Since $K_1$ is a Sylow $2$-subgroup of $K$,
there exists some $w\in W$ such that $L^w\leqslant K_1$.
It follows that $L\unlhd K_1$,
and then $L=1$,
which is impossible. Hence $L\cap K=1$, and so
$L\leqslant K_1H$, and  $L\cong\ZZ_2$.
Thus $\soc(X)=W\times L$,  as in part (c).

By the above  paragraph, we obtain that $\C_{X}(W)=W\times L$.
Let $\ov X=X/L$, and $\ov G=GL/L\cong G$. Let $\ov K_1=K_1L/L\cong K_1$.
Write $\ov G=\ov W{:}\ov H$. Since $H$ normalizes $K_1$,
we conclude that $\ov H$ normalizes $\ov K_1$.
Note that $\ov K_1 \ov H$ acts irreducibly and faithfully on $\ov W$.
By Clifford's Theorem, $\ov W$ can be decomposed as
$\ov W=e(U_1\oplus U_2\oplus\cdots\oplus U_t)$
such that $\ov K_1$ normalises each $U_i$, and all $U_i$ are pairwise non-equivalent and
irreducible with respect to the action of $\ov K_1$.
Recall that $K_1$ is of order at least $4$.
It implies that $t\geqslant 2$.
Let $V_i=eU_i$ for each $i$. Rewrite $\ov W=V_1\oplus V_2\oplus\cdots\oplus V_t$.
Now $\ov H$ normalises $\ov K_1$, we conclude that $\ov H$ preserves such decomposition.
Since the maximal subgroup preserving such decomposition in $\GL(\ov W)$ is $\GL(V_1)\wr\S_t$,
implying that $\ov H$ belongs to $\GL(V_1)\wr\S_t$,
forcing $\ov H$ is imprimitive on $\ov W$.
By \cite[Proposition~2.8]{DF}, we are done, as in part (d).
\qed

We now determine the graph $\Ga$ for the case where $W\unlhd X$, and
$\Ga$ is a normal cover of $\Ga_W$.

\begin{lemma}\label{Le}
Assume that $\Ga$ is a normal cover of $\Ga_W$. Then we have
\begin{itemize}
\item[(i)]
$G$ is normal in $X$, or
\item[(ii)]
$G$ is not normal in $X$, and
\begin{itemize}
\item[(a)]
$\Ga_W\cong\C_{\frac{n}{2}[2]}$, and $n\equiv 0\,(\mod 4)$;

\item[(b)]
$X=W{:}((NH).\calO)$, $X_1\leqslant N.\calO$, $N\cap H\cong\ZZ_2$, and $H$ normalizes $N$, but $H$ does not centralise $N$,
where $N\cong\ZZ_2^l$ with $2\leqslant l\leqslant \frac{n}{2}$, and $\calO\cong1$ or $\ZZ_2$;

\item[(c)]
$W$ is unique and minimal in $X$, and $H$ is imprimitive on $W$;

\item[(d)]
$X/(WN)\cong\ZZ_{\frac{n}{2}}$ or $\D_n$, and $\Ga$ is $X$-arc-transitive if and only if
$X/(WN)\cong\D_n$.
\end{itemize}
\end{itemize}
\end{lemma}
\proof
Let $\ov H=G/W$ with $\ov H=\l\ov h\r$.
Since $\Ga$ is a Cayley graph of $G$,
$\Ga_W$ is a Cayley graph of $\ov H$.
By \cite[Theorem~1.2]{Baik}, either $\Ga_W$ is a normal Cayley graph
or $\Ga_W=\C_{{\frac{n}{2}}[2]}$, and $4$ divides $n$.
It follows that either $G$ is normal in $X$ or $\Ga_W\cong\C_{{\frac{n}{2}}[2]}$.
Suppose that $G$ is not normal in $X$. Then $\Ga_W\cong\C_{{\frac{n}{2}}[2]}$, as in part (a).

Clearly, $\Aut\Ga_W\cong\ZZ_2^{\frac{n}{2}}{:}\D_n$.
Let $\ov K\unlhd \Aut\Ga_W$ such that $\ov K\cong\ZZ_2^{\frac{n}{2}}$.
%Now we choose $\ov K\leqslant\Aut\Ga_W$ such that
%Let $\ov K=\l\s_1,\s_2,\cdots,\s_{{\frac{n}{2}}}\r$.
%Since $\ov H$ is regular on $\Ga_W$, we have $\ov h^{{\frac{n}{2}}}=\s_1\s_2\cdots\s_{{\frac{n}{2}}}$.
Then we may write $\Aut\Ga_W=\ov K~\ov HO$, where $O\cong\ZZ_2$.
Let $B$ be a vertex of $\Ga_W$, and $1\in B$.
Choose $\ov M\leqslant \ov K$ such that $M\cong\ZZ_2^{\frac{n}{2}-1}$ and
$(\Aut\Ga_W)_B=\ov MO$.

Let $\ov X=X/W$.  Since $\ov X~\ov K/\ov K\cong \ov H\calO/(\ov H\calO\cap \ov K)$ where $\calO=1$ or $O$,
we conclude that $\ov X=(\ov X\cap \ov K)\ov H\calO$,
and $\Ga$ is $X$-arc-transitive if and only if $\calO=O$.
Let $\hat K=\ov X\cap \ov K$.
Then $\hat K\unlhd\ov X$, and $\hat K\cap \ov H\cong\ZZ_2$.
Thus $X=W.((\hat K\ov H).\calO)$.
Let $K$ be the preimage of $\hat K$, under $X\to X/W$.
Note that $G$ is a Frobenius group, we conclude that the order of $W$ is odd.
By Hall's Theorem, $K=W{:}N$, where $N\cong \hat K$. It further implies that
$N\cong\ZZ_2^l$, where $l\leqslant \frac{n}{2}$.

Now $(|N|,|W|)=1$, we get $X/W=\N_{X/W}(NW/W)=\N_X(N)W/W$, and so $X=W\N_X(N)$.
Since $H\leqslant X$, it follows that $H^w$ belongs to $\N_X(N)$ for some $w\in W$.
Without loss of generality, we may assume that $H$ belongs to $\N_X(N)$.
Thus $X=W{:}((NH).\calO)$. By comparing the order, we conclude that $N\cap H\cong\ZZ_2$.
If $l=1$, then $NH=H$, and so $G\unlhd X$,
which contradicts the assumption that $G$ is not normal
in $X$.  Then $l\geqslant 2$. Thus $2\leqslant l\leqslant \frac{n}{2}$.

Set $Y=W{:}(NH)$.
Clearly, $G$ is a subgroup of $Y$. Then $Y=GY_1$.
Note that $|Y|=\frac{|W||H||N|}{|H\cap N|}$, and $|Y|=|G||Y_1|$,
we obtain $|Y_1|=\frac{|N|}{|N\cap H|}=\frac{|N|}{2}$.
Let $\ov Y=Y/W$.
Since $Y_1W/W=\ov Y_B$, we have $Y_1W/W\leqslant \ov Y\cap \ov M\leqslant\hat K$,
and hence $Y_1W/W\leqslant NW/W$.
Consequently, $Y_1\leqslant N^{\ov w}$ for some $\ov w\in W$.
For simplicity, we may assume that $Y_1\leqslant N$.
For that case, $Y_1$ has index $2$ in $N$,
and hence $X_1\leqslant N.\calO$,  as in part~(b).

%Let $H=\l h\r$. Then $N=\l Y_1,h^{\frac{n}{2}}\r$.
Let $C{:}=\C_{NH}(W)$. Assume that $C\not=1$.
Clearly, $C$ is normal in $Y$.
Without loss of generality, $C$ is minimal in $Y$.
Since $H$ acts fixed-point-freely on $W$,
we have $C\cap H=1$. % and so $C\cap G=1$.
%Recall that $G\not\cong\D_{2p}$, we conclude that $\Ga$ is a cover of $\Ga_C$.
%It follows that $|C|$ is a divisor of $n$.
Let $\ov C$ be the image of $C$ under $X\to X/W$.
Then $\ov C$ is normal and minimal in $\ov Y$, and hence
$\ov C$ is a subgroup of $\hat K$.
It implies that $C\cong\ZZ_2^{\ell}$ for some $\ell$.

Let $\ov K=\prod_{i=1}^{\frac{n}{2}}\l\s_i\r$.
Note that $\ov H$ acts on $\ov K$ by permuting transitively on all $\s_i$.
Relabeling if necessary, we may assume $\ov h=\s\pi$, where $\s\in\ov K$, and $\pi=(12\cdots\frac{n}{2})^{-1}$.
Let $\ov K_B=\prod_{i\not=1}\l\s_i\r$.
Choose $\ov B,\tilde B\in\Ga(B)$ such that $\ov K_{\ov B}=\prod_{i\not=2}\l\s_i\r$
and $\ov K_{\tilde B}=\prod_{i\not=\frac{n}{2}}\l\s_i\r$.
Pick some $x\in C$ such that $x=\s_{i_1}\cdots\s_{i_k}$, where $i_1=2$, and
$2\leqslant i_j\leqslant i_{j+1}\leqslant\frac{n}{2}$.
Then $\l x,x^{\ov h^{\frac{n}{2}-i_k}}\r\leqslant \ov K_B\cap C$.
It follows that $\Ga_W$ is $\ov C{:}\ov H$-edge-transitive Cayley graph of valency $4$.

Let $Z=(W\times C){:}H$. By the above paragraph, $\Ga$ is $Z$-edge-transitive.
However, $\Ga$ is not $Z$-arc-transitive.
By Lemma~\ref{Cos}, $\Ga=\Cos(Z,Z_1,Z_1\{g,g^{-1}\}Z_1)$.
It is obvious that $Z_1\cong\ZZ_2^\ell$.
Now we may assume $Z_1=\l\t_1,\t_2,\ldots,\t_{\ell}\r$.
Write $h=\l(h_1,h_2)\r$.
If $\t_i=(1,c)$ for some $i$, since $C$ is minimal in $Z$,
we conclude that $Z_1=C$, which is impossible.
Thus each $\t_i$ has the form $(1,c)h^{\frac{n}{2}}$ or $(u,c)h^{\frac{n}{2}}$, where $u\in W\setminus\{1\}$,
and $c\in C$. Note that $n$ is divisible by $4$.
Then we may write $g$ as $(v,c')h$, where $v\in W$, and $c'\in C$.
Assume $\t_{i_0}=(u_0,c)h^{\frac{n}{2}}$ for some $i_0$, where $u_0\not=1$.
Then all $\t_i$ have the form $(u_0,c')h^{\frac{n}{2}}$, where $c'\in C$.
If $\t_i^g=\t_j$ for some two $i,j$, it is easy to show that $\l Z_1,g\r$ is a subgroup of $C{:}\l g\r$,
which is a contradiction.
Thus all $\t_i^{g}$ do not belong to $Z_1$,
which leads to the valency of $\Ga$ is greater than $4$, again
a contradiction.
Similarly, we also exclude the other case.
Thus $C=1$, namely, $NH$ acts faithfully and irreducibly on $W$.
That is to say, $W$ is the unique minimal normal subgroup of $X$.
Arguing similarly as Lemma~\ref{Lee}, we obtain that
$H$ is imprimitive on $W$, as in part~(c).

Clearly, $\Ga_{WN}$ is a cycle.  Since $X/(WN)$ is transitive on $V\Ga_{WN}$,
we conclude that $X/(WN)\cong\ZZ_{\frac{n}{2}}$ or $\D_n$.
For that case, $\Ga$ is arc-transitive if and only if $X/(WN)\cong\D_n$, as in part~(d).
\qed

With the above preparation, we are ready to embark on the proof of Theorem~\ref{soluble}.

\vskip0.1in
{\bf Proof of Theorem~\ref{soluble}:}
If $G\lhd X$, then by Lemma~\ref{G}, we have  $X_1\leqslant \D_8$, as in Theorem~\ref{soluble}~(1).
In what follows, we assume that $G$ is not normal in $X$.

Assume first that $p$ is an odd prime.
By Lemmas~\ref{D_2p}-\ref{Le}, if  $W$ is not normal in $X$,
we obtain that  $\Ga\cong\C_{p[2]}$, and  $\Aut\Ga\cong\ZZ_2^p{:}\D_{2p}$, as in Theorem~\ref{soluble}~(2).
If $W$ is normal in $X$, and $\Ga_W$ is a cycle, by Lemma~\ref{Lee}, part~(3) of Theorem~\ref{soluble} occurs.
If $W$ is normal in $X$, and $\Ga$ is a cover of $\Ga_W$,
from Lemma~\ref{Le}, it follows that part~(4) of Theorem~\ref{soluble} holds.

Assume now that $p=2$. %By Lemma~$\ref{Lei}$, $F=O_2(X)$, and $W\leqslant F$.
By Lemmas~\ref{Lei} and \ref{Lee}, Theorem~\ref{soluble}\,(5) occurs.
\qed

\section{Insoluble automorphism groups}
Let $G=W{:}H\cong\ZZ_p^d{:}\ZZ_n$ be a primitive Frobenius group.
Assume that $\Ga=(V\Ga,E\Ga)$ is  a connected
$X$-edge-transitive tetravalent Cayley graph of $G$,
where $G\leqslant X\leqslant\Aut\Ga$.
In this section, we study the case where the automorphism group $X$ is insoluble.

For a finite group $R$, the socle of $R$, denoted by $\soc(R)$, is the subgroup generated by all minimal normal subgroups of $R$.
The group $R$ is said to be almost simple if its socle $\soc(R)$ is a non-abelian simple group.

\vskip0.1in
\centerline{\bf TABLE~2:\ {\rm Almost simple automorphism groups.}}

\[\begin{array}{|c|c|c|c|c|}\hline
X& G & X_1   \\ \hline
\PSL(3,3){:}\ZZ_2 &\D_{26}  & \ZZ_3^2{:}\GL(2,3)\\ \hline
\PGL(2,7) & \D_{14} & \S_4 \\\hline
\PGL(2,7) & \ZZ_7{:}\ZZ_3 & \D_{16} \\\hline
\PGL(2,7) & \ZZ_7{:}\ZZ_6 & \D_8  \\\hline
\PSL(2,23) &\ZZ_{23}{:}\ZZ_{11}  & \S_4\\ \hline
\PSL(2,11) &\ZZ_{11}{:}\ZZ_{5}  & \A_4 \\ \hline
\PGL(2,11) &\ZZ_{11}{:}\ZZ_{5}  & \S_4 \\ \hline
\PGL(2,11) &\ZZ_{11}{:}\ZZ_{10}  & \A_4\\ \hline
\end{array}\]

\vskip0.1in
We now determine the structure of insoluble group $X$.
Denote by $R(X)$ the maximal solvable normal subgroup of $X$.
We first treat the case where $R(X)=1$.

\begin{lemma}\label{Z(X)}
Let $N$ be minimal and normal
in $X$. If $R(X)=1$, then $\C_X(N)=1$.
\end{lemma}
\proof Note that $N$ is minimal in $X$. Since $R(X)=1$,
we have $N\cong T^k$, where $T$ is a nonabelian simple group,
and $k$ is an integer.
Clearly, $\Z(N)=1$.  Let $C{:}=\C_X(N)$.
Since $N\unlhd X$, we have $C\unlhd X$. Suppose that $C\not=1$.
By our assumption, $C$ is insoluble.
Notice that $N\cap G\unlhd G$, we conclude that $N\cap G=1$ or $W\leqslant N\cap G$.
For the former, $|N|$ divides $|X_1|$, and so $N$ is soluble, contrary to our assumption.
Thus $W\leqslant N\cap G$.  Similarly,  $W\leqslant C\cap G$.
It follows that $W\leqslant N\cap C$, a contradiction. Thus $C=1$.
\qed

\begin{lemma}\label{X}
If $R(X)=1$, then $X$ is almost simple.
\end{lemma}
\proof By Frattini argument, we have that $X=GX_u$, where $u\in V\Ga$.
By Lemma~$\ref{Lu}$, either $X_u$ is a $2$-group
or  $|X_u|$ divides $2^43^6$.
Let $N$ be a minimal normal subgroup of $X$.
By our assumption, $N$ is unsolvable.
So $N=T_1\times T_2\times\cdots\times T_k$, where $T_i\cong T$ is a nonabelian simple group
for any $i$.
By \cite{L.Kazarin}, we obtain that $T$ is one of the following:
\[\PSL(2,q)(q>3),\, \PSL(3,q)(q<9),\,\PSL(4,2),\, \PSp(4,3),\ \PSU(3,8),\ \mbox{or}\ \M_{11}.\]

In what follows, suppose that $k\geqslant2$.
Since $W$ is minimal and normal in $G$ and $N\cap G\not=1$, we conclude that $W\leqslant N$.
Let $r>3$ be a prime divisor of $|T|$.
Since $r$ divides $|X|$ and $(|W|,|H|)=1$,
we conclude that $r$ divides either $|W|$ or $|H|$.
Suppose first that $r$ divides $|W|$.
Then $T_i\cap W\not=1$ for each $i$.
Let $W_i=T_i\cap W$ for $1\leqslant i\leqslant k$.
Assume that $N\cap H=1$.
Then $\frac{|N|}{|W|}$ divides $|X_u|$.
So does $\prod_{i=1}^k\frac{|T_i|}{|W_i|}$.
An inspection of the above simple groups shows $T=\A_5$ and $k=2$.
By Lemma~\ref{Z(X)}, $X\lesssim\S_5\wr\ZZ_2$.
For this case, $\frac{|N|}{|W|}=144$.
By Lemma~\ref{Lu}, the only possibility is that $X\cong\A_5\times\A_5$.
Clearly, this is a contradiction.

Thus $N\cap H\not=1$. Let $\ov H=N\cap H$.
Since $G$ is a Frobenius group, it follows that $\ov H$ is a diagonal subgroup of $N$.
Write $\ov H=\l\s_1\s_2\cdots\s_k\r$ where $\l \s_i\r\cong\ov H$, and $\s_i\in T_i$ for each $i$.
Let $H_i=\l\s_i\r$ where $1\leqslant i\leqslant k$.
Then $G_i=W_i{:}H_i$  is a Frobenius group.
For this case, we obtain that $\frac{|N|}{|W||\ov H|}$ divides $|X_u|$.

Let $T=\PSL(2,q)$ where $q>3$. By \cite[Theorem~6.25]{Suzuki},
we conclude that $ G_i\leqslant [q]{:}[\frac{q-1}{d}]$,
$G_i\leqslant \D_{\frac{2(q\pm1)}{d}}$, $G_i\leqslant\A_5$, or $G_i\leqslant\PGL(2,r)$,
where $d=(2,q-1)$, and $r\div q$.
Suppose $ G_i\leqslant [q]{:}[\frac{q-1}{d}]$. Then $W{:}\ov H\leqslant [q]^k{:}[\frac{q-1}{d}]$,
namely, $N\cap G\leqslant [q]^k{:}[\frac{q-1}{d}]$.
Since $|N{:}N\cap G|$ divides $|X_u|$, we conclude that $\frac{d|N|}{q^k(q-1)}$ divides $|X_u|$.
If $q$ is even, then $(q+1)^k(q-1)^{k-1}\div2^43^6$,
which is a contradiction because $q+1$ and $q-1$ are two distinct odd numbers.
If $q$ is odd, then $(q+1)^k(\frac{q-1}{2})^{k-1}\div 2^43^6$.
By easy calculations, $q=5$ and $k=2$.
For this case, the only possibility is that $G\cong\ZZ_5^2{:}\ZZ_8$ and $X\cong\S_5\wr\ZZ_2$.
By Lemma~\ref{Lu}, we have that $X_u\cong\S_3\times\S_4$.
Without loss of generality, we may assume that $X=(\S_5\times\S_5){:}\l \s\r$, where
$\s$ permutes the first and second coordinates.
Note that $G\cap X_u=1$.
By MAGMA \cite{Magma}, there is an element $G$ (up to conjugate) in $X$, and
there are two elements $X_u$ (up to conjugate) in $X$ such that their intersections equal to 1.
Choose $w,h\in \S_5$ such that $o(w)=5$, $o(h)=4$ and $w^h=w^2$.
For this case, write $G=W{:}H$ with $W=\l(w,1),(1,w)\r$ and $H=(h,1)\s$.
Meanwhile, we choose $X_u=\l((123),1),((12),1),(1,(1234),(1,(12)))$ or $\l((123),1),((12)(45),1),(1,(1234)),(1,(12))\r$.
It is simple to show that, for the above two choices, $X_u$ belongs to different conjugate classes of $X$, and
$X_u\cap G=1$.
%By MAGMA \cite{Magma}, there are two conjugate classes in $X$,
%and  for each element $L$ of these two conjugate classes, $L$ is isomorphic to $\S_3\times\S_4$  and $L\cap G=1$.
Choose $v\in\Ga(u)$.
By Lemma~\ref{Cos}, write $\Ga=\Cos(X,X_u,X_uoX_u)$, where $o\in\N_X(X_{uv})\setminus X_{u}$ and $o^2\in X_{uv}$.
Since $\Ga$ is $X$-arc-transitive graph, we conclude $|X_u{:}X_{uv}|=4$, and hence $|X_{uv}|=36$.
In such two cases,
again by MAGMA \cite{Magma}, there is no $o\in\N_X(X_{uv})$ such that
$\l X_u,o\r=X$, namely, $\Ga$ is not connected.
Suppose that $G_i\leqslant\D_{\frac{2(q\pm1)}{d}}$ or $G_i\leqslant\A_5$
where $1\leqslant i\leqslant k$.
Arguing similarly as above, we conclude that $q=4$, $k=2$,
and $G_i\cong\D_{10}$ for each $i$.
For this case, the only possibility is that $G\cong\ZZ_5^2{:}\ZZ_8$ and $X\cong\S_5\wr\ZZ_2$,
which is impossible by the above discussion.
Thus $G_i\leqslant \PGL(2,r)$. Write $q=p^n$ where $p$ is a prime, and $n$ is a natural number.
Then $r=p^m$, where $m\div n$. Let $n=ms$ with $s\geqslant2$.
Note that $|N{:}N\cap G|$ divides $|X_u|$. So does $|T{:}\PGL(2,r)|^k$.
It follows that $\frac{1}{2}p^{m(s-1)}(\Sigma_{i=1}^sp^{m(s-i)})(\Sigma_{j=1}^sp^{2m(s-j)})$ divides $2^43^6$.
By easy calculations, there are no $p,n$ and $m$ satisfying the above equation.
Thus this case is excluded.

Let $T=\PSL(3,q)$ with $q<9$. Assume that $q=2$. By Atlas~\cite{atlas}, we have $G_i\cong\ZZ_7{:}\ZZ_3$  where $1\leqslant i\leqslant k$.
Clearly, $\frac{|N|}{|W||\ov H|}$ does not divide $|X_u|$. Similarly, we can exclude the other cases.
Let $T=\PSL(4,2)$. Again by Atlas~\cite{atlas}, we conclude that $35\notdiv|G_i|$ where $1\leqslant i\leqslant k$.
It implies that $5$ or $7$ divides $|X_u|$, which is impossible.
Arguing similarly as above, we can conclude $T$ can not equal to other simple groups.

Suppose now that $r$ divides $|H|$. If $r\notdiv |H_i|$, then $r$ divides $|X_u|$, which is impossible.
Thus $r$ divides $ |H_i|$ for each $i$.  Recall that $\frac{|N|}{|W||\ov H|}$ divides $|X_u|$,
we conclude that $r$ divides $|X_u|$, again a contradiction.
Therefore, $X$ is almost simple.
\qed

\begin{lemma}\label{ll}
Let $X$ be an almost simple group with $\soc(X)=\PSL(2,7)$.
Assume that $\Ga$ is not $(X,2)$-arc-transitive.
Then either $X=\PGL(2,7)$,  $X_1=\D_8$ and $G=\ZZ_7{:}\ZZ_6$, or
$X=\PGL(2,7)$, $X_1=\D_{16}$ and $G=\ZZ_7{:}\ZZ_3$.
\end{lemma}
\proof
By Frattini argument, we have that $X=GX_u$, where $u\in V\Ga$.
Since $\Ga$ is  not $(X,2)$-arc-transitive, it follows that $X_u$ is a $2$-group.
Note that $G$ is a Frobenius group.  Checking
the subgroups of $\PGL(2,7)$ in the Atlas \cite{atlas}, we obtain $G=\ZZ_7{:}\ZZ_6$, or $\ZZ_7{:}\ZZ_3$.

Denote by $T$ the socle $\soc(X)$.
Assume first that $G=\ZZ_7{:}\ZZ_6$. Since $\ZZ_7{:}\ZZ_3$ is maximal in $T$,
we have $X=\PGL(2,7)$.
It follows that $X_u\cong\D_8$.
Assume now that $G=\ZZ_7{:}\ZZ_3$. Furthermore, assume that $X=\PSL(2,7)$.
Then $\Ga$ is a connected tetravalent $X$-edge-transitive Cayley graph,
and $X_u\cong\D_8$ is a Sylow $2$-subgroup
of $X$. Choose $v\in\Ga(u)$. Then $|X_u{:}X_{uv}|=2$ or $4$.
Since $\Ga$ is vertex transitive graph, we write
$\Ga$ as coset graph $\Cos(X,H,H\{x,x^{-1}\}H)$, where $H=X_u\cong\D_8$, and $x\in X$ is such that $\langle H, x\rangle=X$;
in particular, $x\notin H$.

Suppose that $|X_u{:}X_{uv}|=4$.  Then $\Ga$ is $X$-arc-transitive.
By Lemma~\ref{Cos}, we
choose $x$ such that $(u,v)^x=(v,u)$, resulting $x\in \N_X(X_{uv})\cong\D_8$. In particular, $\N_X(X_{uv})\not=X_u$.
Then $|\N_{X_u}(X_{uv})|=4$. Hence $\N_{X_u}(X_{uv})$ is normal in both $\N_{X}(X_{uv})$ and $X_u$,
and so $\N_{X_u}(X_{uv})\unlhd\langle X_u,\N_{X}(X_{uv})\rangle$.
Checking the subgroups of $\PSL(2,7)$,
we obtain that $\langle X_u,\N_{X}(X_{uv})\rangle\cong\S_4$,
which contradicts the fact $\langle X_u, x\rangle=X$.

Suppose that $|X_u{:}X_{uv}|=2$.
Then $|X_{uv}|=4$, and hence $X_{uv}\unlhd M=
\langle X_u, X_v\rangle$, so $M\cong\S_4$.
 By Lemma~\ref{Cos},
we may choose $x$ such that $u^x=v$.
It is clear that $X_u$
and $X_v$ are two Sylow $2$-subgroup of $M$,
there exists some $y\in M$ such that $X_u^y=X_v=X_u^x$. Hence $xy^{-1}\in\N_X(X_u)=X_u$,
so $\langle X_u, x\rangle\leqslant\langle X_u, xy^{-1}, y\rangle\leqslant M$,
again a contradiction.
Thus $X=\PGL(2,7)$.
\qed

Lemma~\ref{X} tells us that if $X$ is insoluble and $R(X)=1$, then $X$ is almost simple.
The next two lemmas determine $\Ga$ for the case where $X$ is almost simple.

\begin{lemma}\label{al}
If $X$ is almost simple, then,  for $u\in V\Ga$,  the triple $(X, G, X_u)$ is one of
the triples listed in Table~$2$.
\end{lemma}
\proof
By Frattini argument, we have that $X=GX_u$.
Since $\Ga$ is a connected Cayley graph of valency $4$,
we obtain that $X_u$ is a $\{2,3\}$-group.
It follows that $X$ is decomposed as a product of two solvable subgroups.
Denote by $T$ the socle $\soc(X)$.
By \cite{L.Kazarin}, we conclude that $T$ appears in Lemma~$\ref{X}$.
In what follows, we process our analysis by two cases.

\vskip0.07in
{\bf Case 1:}\ Assume that $\Ga$ is $(X,2)$-arc-transitive.
By Lemma~$\ref{Lu}$, $|X_u|$ divides $2^43^6$.

Assume first that $T=\PSL(2,q)$ with $q>3$.
%Then $X=\PSL(2,q).\calO$, where $\calO\leqslant\Out(T)$.
If $q=5$, then $G=\ZZ_5$ and $X_u=\A_4$ or $\S_4$, a contradiction.
If $q=7$, by the Atlas \cite{atlas}, we conclude that $X=\PGL(2,7)$, $G=\D_{14}$, and $X_u=\S_4$.
Suppose  $q=11$.  Checking the subgroups of $\PGL(2,11)$ in the Atlas~\cite{atlas},
we know that $X=\PSL(2,11).\calO$, $G=\ZZ_{11}{:}(\ZZ_5\times\calO_1)$, and $X_u=\A_4.\calO_2$,
where $\calO=\ZZ_2$, and $\calO_1\calO_2=\calO$, (refer to \cite[Theorem~1.5]{Maru}).
Suppose that $q=23$. By \cite[Theorem~1.1]{Zai Ping Lu}, $X=\PSL(2,23)$, $G=\ZZ_{23}{:}\ZZ_{11}$ and $X_u=\S_4$.

In what follows, we assume that $q\not=4,5,7,11,23$.
Then, by \cite[Proposition~4.1]{Xia}, interchanging $G$ and $X_u$ if necessary,
$G\cap T\leqslant\D_{2(q+1)/d}$ and $[q]\unlhd T\cap X_u\leqslant [q]{:}[\frac{q-1}{d}]$
where $d=(2,q-1)$.
Let $T_u=T\cap X_u$.
Assume that $G\cap T\leqslant\D_{2(q+1)/d}$, and $[q]\unlhd T_u\leqslant  [q]{:}[\frac{q-1}{d}]$.
Then $\frac{q(q-1)}{2}$ divides $|T{:}T\cap G|$.
Since $|T{:}T\cap G|=|TG{:}G|$ divides $|X_u|$, we conclude that $\frac{q(q-1)}{2}\div 2^43^6$.
By easy calculations, $q=9$.
It follows that $T_u^{\Ga(u)}\cong\ZZ_3$ or $\S_3$.
However, $T_u^{\Ga(u)}\unlhd X_u^{\Ga(u)}\leqslant \S_4$, and so $X_u^{\Ga(u)}\cong \S_3$,
which is a contradiction because $X_u^{\Ga(u)}$ is $2$-transitive on $\Ga(u)$.
Thus $[q]\leqslant T\cap G\leqslant  [q]{:}[\frac{q-1}{d}]$ and $T_u\leqslant \D_{2(q+1)/d}$.
Then  $T_u^{\Ga(u)}\cong\S_3$ or $\ZZ_3$, and so $ X_u^{\Ga(u)}\cong \S_3$,
again a contradiction.

Assume now that $T=\PSL(3,q)$ ($q<9$), $\PSp(4,3)$, $\PSL(4,2)$, $\PSU(3,8)$  or $\M_{11}$.
It is clear that $T\cap G\not=1$. Then $W\leqslant T$.
Suppose that $T=\PSp(4,3)$.
Using \cite[Proposition~4.1]{Xia},
$T\cap G=\ZZ_2^4{:}\ZZ_5$ and $T\cap X_u=3_{+}^{1+2}{:}2\A_4$.
For that case, it is clear that $G\cap X_u\not=1$, a contradiction.
Suppose that $T=\M_{11}$. Then $X=\M_{11}$. Again by \cite[Proposition~4.1]{Xia},
$G=\ZZ_{11}{:}\ZZ_5$, and $X_u=\M_9.2$, again a contradiction, refer to Lemma~\ref{Lu}.
Suppose that $T=\PSL(3,4)$. By Atlas~\cite{atlas}, we conclude that $35\notdiv |T\cap G|$.
It implies that $5$ or $7$ divides $|X_u|$, which is impossible.
Similarly, $T$ does not equal $\PSL(3,q)$ with $5\leqslant q\leqslant 8$,  $\PSL(4,2)$, or $\PSU(3,8)$.
If $T=\PSL(3,3)$,  we obtain that $X=\PSL(3,3){:}\ZZ_2$, $G=\D_{26}$ and
$X_u=\ZZ_3^2{:}\GL(2,3)$, refer to \cite[Proposition~4.1]{Xia}.

\vskip0.07in
{\bf Case 2:}\ Assume that $\Ga$ is not $(X,2)$-arc-transitive.
By Case 1, we only need deal with the case where $T=\PSL(2,q)$ with $q>3$.
Since $\Ga$ is not $(X,2)$-arc-transitive, $X_u$ is a 2-group. Since $X=GX_u$ and $G\cap X_u=1$,
$G$ contains a $2'$-Hall subgroup of $X$. Then $G\cap T$ contains a $2'$-Hall subgroup of $T$.
By Lemma~\ref{Kazarin}, we have that $T=\PSL(2,q)$,
$T\cap G=\ZZ_q{:}\ZZ_{\frac{q-1}{2}}$, and $T_u=\D_{q+1}$,
where $q=2^e-1$ is a prime (see \cite[Proposition~4.1]{Xia}). Suppose that $e=3$. Then $q=7$.
By Lemma~\ref{ll}, the statement follows.
In what follows, we assume that $e\geqslant5$.

Note that $\ZZ_q{:}\ZZ_{\frac{q-1}{2}}$ is maximal in $T$.
By \cite[Theorem~1.1]{Hua Zhang}, we conclude that
$G=\ZZ_q{:}\ZZ_{q-1}$, and hence $X=\PGL(2,q)$, and $X_u=T_u=\D_{q+1}$.
%By the above paragraphs, we have $|T{:}(T\cap G)|=|X_u|=q+1$,
%and so $T_u=X_u=\D_{q+1}=\D_{2^e}$.
Let $v\in\Ga(u)$. By Lemma~\ref{Cos}, $T_{uv}$ has index $2$ or $4$ in both $T_u$ and $T_v$.
Since $e\geqslant5$, $T_{uv}$ contains a subgroup $C\cong\ZZ_4$. It is easily shown that $C$ is normal in both $T_u$
and $T_v$, and so $C\lhd L{:}=\l T_u,T_v\r$.
In view of \cite[p.417]{Suzuki}, each Sylow $2$-subgroup of $T$ is maximal in $T$.
Note that $T_u$ is a Sylow $2$-subgroup of $T$,
it follows that $L=T_u=T_v$.
By the connectedness of $\Ga$, we obtain $L$ fixes each vertex of $\Ga$,
which is impossible.
This completes the proof of Lemma~\ref{al}.
\qed

\vskip0.07in

By  \cite[Theorem~1.1]{Zai Ping Lu}, we have the following lemma.

\begin{lemma}\label{T}
Let $X$ be in Table~$2$, and let $T=\soc(X)$.
Then we have:
\begin{itemize}
\item[(1)]
if $T=\PSL(3,3)$, then $\Ga$
is isomorphic to the graph given in Example $\ref{TA}$;

\item[(2)]
if $T=\PSL(2,7)$, then $\Ga$ is isomorphic to a graph given in
Examples~$\ref{TA}$, $\ref{A}$ and $\ref{B}$;

\item[(3)]

if $T=\PSL(2,23)$, then  $\Ga$ is isomorphic to the graph given in
Example~$\ref{C}$;

\item[(4)]

if $T=\PSL(2,11)$,  then $\Ga$  is isomorphic to a graph given in Examples~$\ref{C}$-$\ref{D}$.

\end{itemize}
\end{lemma}

We now begin with treating the case where  $R(X)\not=1$.

\begin{lemma}\label{R}
If $R(X)\cap G=1$, then $G=\ZZ_{11}{:}\ZZ_{10}$ and $X=\PGL(2, 11)\times\ZZ_2$.
\end{lemma}
\proof Let $\ov X_1=X_1R(X)/R(X)$, $\ov G=GR(X)/R(X)$, and $\ov X=X/R(X)$.
Since $X=GX_1$, we conclude that $\ov X=\ov G~\ov X_1$.
Since $R(X)\cap G=1$, implying that $\ov G\cong G$ is a Frobenius group.
Since $\ov X$ is insoluble, $\Ga$ is a cover of $\Ga_{R(X)}$, refer to Lemma~\ref{insoluble}.
Let $B$ be a vertex of $\Ga_{R(X)}$, where $1\in B$.
By Frattini argument, $\ov X=\ov G~\ov X_B$.
Clearly, $\ov X_1\leqslant \ov X_B$, and $\ov G\cap\ov X_B\not=1$.

Assume that $\ov X$ is not almost simple. Let $\ov N$ be a minimal normal subgroup of $\ov X$.
Arguing similarly as Lemma~\ref{Z(X)}, we obtain that $\C_{\ov X}(\ov N)=1$.
Suppose that $\soc(\ov X)=\A_5\times\A_5$. By the definition of $G$, we conclude that $\ov G\cong\ZZ_5^2{:}\ZZ_8$.
For this case, $\ov X=((\A_5\times\A_5){:}\ZZ_2){:}\ZZ_2$ is a subgroup of $\S_5\wr\S_2$, and $\ov X_B\cong\S_3\times\S_4$.
It implies that $\ov G\cap\ov X_B\cong\ZZ_2$.
%Note that $|R(X)|$ divides $|G|$, we conclude that $|R(X)|$ is a $2$-power.
%Since $|X|=|R(X)||\ov X|$ and $|\ov X|=\frac{|\ov G||\ov X_B|}{2}$, we conclude that $|X_1|=|\ov X_B|$.
By MAGMA \cite{Magma}, there are two elements $\ov G$ (up to conjugate) in $\ov X$,
and there is an element $\ov X_B$ (up to conjugate)  in $\ov X$
such that their intersections are isomorphic to $\ZZ_2$.
Choose $\ov B\in\Ga(B)$.
By Lemma~\ref{Cos}, write $\Ga_{R(X)}=\Cos(\ov X,\ov X_B,\ov X_Bo\ov X_B)$,
where $o\in\N_{\ov X}(\ov X_{B\ov B})\setminus \ov X_{B}$ and $o^2\in \ov X_{B\ov B}$.
Since $\Ga$ is $\ov X$-arc-transitive graph, we conclude $|\ov X_B{:}\ov X_{B\ov B}|=4$, and hence $|X_{B\ov B}|=36$.
Again by MAGMA \cite{Magma}, for each choice of $\ov G$ and $\ov X_B$,
there is no $o\in\N_{\ov X}(\ov X_{B\ov B})$ such that
$\l \ov X_B,o\r=\ov X$, namely, $\Ga_{R(X)}$ is not connected.
For other cases, arguing similarly as Lemma~\ref{X},
we exclude these possibilities. Thus $\ov X$ is almost simple.

Let $\ov T=\soc(\ov X)$. By \cite{L.Kazarin}, we obtain that $\ov T$ is one of the following:
\[\PSL(2,q)(q>3),\, \PSL(3,q)(q<9),\,\PSL(4,2),\, \PSp(4,3),\ \PSU(3,8),\ \mbox{or}\ \M_{11}.\]
Suppose first that $\ov T=\PSL(2,q)$ where $q=5,7,11,23$.
If $q=5$, the only possibility is that $\ov G\cong\ZZ_5{:}\ZZ_4$, refer to \cite[Theorem~1.1]{Zai Ping Lu}.
For this case, $\Ga$ is a Cayley graph of order $20$,
by \cite[Theorem~5.3]{Pan}, $G$ is normal in $X$, which is a contradiction.
%For the latter, $|G|$  is square-free, which is impossible, refer to \cite[Theorem~1.1]{Zai Ping Lu}.
If $q=7,11$ or $23$, then $|G|$ is square-free, refer to  Atlas~\cite{atlas}.
%By Atlas~\cite{atlas}, we conclude that $\ov X=\PGL(2,11)$, $\ov G=\ZZ_{11}{:}\ZZ_{10}$, and $\ov X_1=\A_4$.
Again by \cite[Theorem~1.1]{Zai Ping Lu},
 $X=\PGL(2,11)\times\ZZ_2$, $G=\ZZ_{11}{:}\ZZ_{10}$, $X_1=\S_4$
and $\Ga$  is isomorphic to the graph given in Example~\ref{D}.
Arguing similarly as Lemma~\ref{al} with $\ov X=\ov G~\ov X_B$ in the place $X=GX_1$,
we can exclude other cases.
This completes the proof.
\qed

\begin{lemma}\label{Rr}
If $R(X)\cap G\not=1$, then we have:
\begin{itemize}
\item[(a)]
$G\cong\ZZ_p^4{:}\ZZ_5$, $X=W.\overline X$, and $\Ga_W=\K_5$, where $\soc(\overline X)\cong \A_5$;
\item[(b)]
$G\cong\ZZ_p^4{:}\ZZ_{10}$, $X=W.(\overline X\times\ZZ_2)$,
and $\Ga_W=\K_{5,5}-5\K_2$, where $\soc(\overline X)\cong\A_5$.
\end{itemize}
\end{lemma}
\proof Let $R=R(X)\cap G$.
Assume that $R(X)\cap G\not=1$.
From the minimality of $W$ in $G$, it follows that $R\geqslant W$.
Since $X/R(X)$ is insoluble, by Lemma~$\ref{insoluble}$, $\Ga$ is a normal cover of $\Ga_{R(X)}$.
Hence $GR(X)/R(X)\leqslant\Aut\Ga_{R(X)}$.

We claim that $R=R(X)$. Let $\ov H=HR(X)/R(X)$.
Since $\Ga$ is a Cayley graph of $G$,
implying that $\Ga_{R(X)}$ is a connected Cayley graph of $\ov H$.
It follows that $|R(X)||\overline H|=|G|$, and hence $|R(X)|=|W||R(X)\cap H|$.
So $R(X)\leqslant G$, and then $R=R(X)$, as claimed.
Furthermore, it implies that $W$ is a normal subgroup of $X$.
By \cite[Theorem~1.2]{Baik}, we obtain that either
$\Ga_W\cong\K_5$
and $H\cong\ZZ_5$ or $\Ga_W\cong\K_{5,5}-5\K_2$ and  $H\cong\ZZ_{10}$.
In the former case, we easily obtain $\Aut\Ga_W\cong\S_5$, and in the latter case, $\Aut\Ga_W\cong\S_5\times\ZZ_2$.

Let $\ov X=X/W$.
Since $\Ga$ is a cover of $\Ga_W$, it follows that
$\ov X$ is a subgroup of $\Aut\Ga_W$.
Assume first that $\Ga_W\cong\K_5$.
Notice that $\ov X$ is insoluble, we conclude that $\soc(\ov X)\cong\A_5$.
Assume now that $\Ga_W\cong \K_{5,5}-5\K_2$.
Since $H\cong\ZZ_{10}$, and $X$ is insoluble, we obtain
$\ov X\cong L\times\ZZ_2$ where $\soc(L)\cong\A_5$.

Recall that $H$ is irreducible on $W$, we conclude that $d=4$, refer to \cite[Lemma~3.1]{Wenqin}.
Hence $G\cong\ZZ_p^4{:}\ZZ_5$ or $\ZZ_p^4{:}\ZZ_{10}$.
\qed

\vskip0.1in
The assertion of Theorem~\ref{isoluble} follows from
Lemmas~\ref{X}-\ref{Rr}.

\section{Half-transitive graphs}
In the last section, we aim to prove Theorem~\ref{solubles}.

Let $p$ be an odd prime, and $d>1$ be an odd integer. Let $n$ be a primitive
divisor of $p^d-1$, and let
\[G=W{:}\langle h\rangle=\ZZ_p^d{:}\ZZ_n<\AGL(1,p^d).\]

\begin{construction}\label{construction}
{\rm
Let $i$ be coprime to $n$ such that $1\leqslant i\leqslant n-1$, and let $a\in W\backslash\{1\}$. Let
\[ \left\{
\begin{array}{ll}
 S_i = \{ah^i, a^{-1}h^i, (ah^i)^{-1}, (a^{-1}h^i)^{-1}\},\\
\Ga_i = \Cay(G,S_i).
\end{array}
\right.\]
}
\end{construction}

\vskip0.1in
{\bf Proof of Theorem~\ref{solubles}:}
Let $X=\Aut\Ga$.
Let $\Ga=\Cay(G,S)$ be connected, edge-transitive and of valency $4$.
Note that $\l h\r$ is primitive on $W$, $d>1$ is odd, and $p$ is an odd prime.
By Theorem~\ref{soluble} and Theorem~\ref{isoluble},
we obtain that  $G$ is normal in $X$.
In view of \cite[Lemma~2.1]{Godsil},
we have $X=G{:}\Aut(G,S)$.

By Lemma~$\ref{G}$, we have that $X_1=\Aut(G,S)\leqslant\D_8$.
By \cite[Proposition~12.10]{Doerk},
$\Aut(G)\cong\ZZ_p^d{:}(\ZZ_{p^d-1}{:}\ZZ_d)$.
Since $d$ and $p$ are odd,  $\Aut(G)$
has a cyclic Sylow $2$-subgroup.
It follows that $X_1=\l\s\r\cong\ZZ_4$ or $\ZZ_2$.
Thus $\s$ fixes an element of $G$ of order $n$, say $f\in G$ such that $o(f)=n$ and
$f^\s=f$. Then $G=W{:}\l f\r$, and $X=G{:}\l\sigma\r=(W{:}\l f\r){:}\l\s\r$.
Moreover, it implies that all involutions of $\Aut(G)$ are conjugate.
Recall that $G$ is a Frobenius group,  every involution of $\Aut(G)$ inverts all elements of $W$.

Since $\Ga$ is connected, $\l S\r=G$ and $\Aut(G,S)$ is faithful on $S$.
Assume that $S$ contains an involution.
Recall that $\Ga$ is $X$-edge-transitive, we conclude that $S$ consists of involutions.
By the proof of Lemma~\ref{G}, $G\cong\D_{2p}$, against our assumption.
Hence $S$ does not contain an involution. For that case,
we may write $S=\{x,x^{-1},y,y^{-1}\}$
such that  either $o(\s)=2$ and $(x,y)^\s=(y,x)$, or  $o(\s)=4$ and $(x,y)^\s=(y,x^{-1})$, refer to~\cite[Proposition~1]{Praeger}.
Now $x=af^i$ , where $a\in W$ and $i$ is an integer. Suppose that $o(\s)=4$.
Then $y=x^{\s}=(af^i)^{\s}=a^\s f^i$, and  $a'f^{-i}=f^{-i}a^{-1}=(af^i)^{-1}=x^{-1}=x^{\s^2}=a^{\s^2}f^i=a^{-1}f^i$.
It follows that $f^{2i}=1$, and hence $f^i$ has order $1$ or $2$.
If $f^i=1$, then $x=a$, and $y=x^\s=a^\s$, belonging to $W$, and so $\l S\r\leqslant W< G$,
which is a contradiction. Thus $f^i$ has order $2$.
Note that $f^i$ inverts each element of $W$, we conclude that $x$ has order $2$,
again a contradiction. Thus $\s$ is an involution,
and so $(x,y)^\s=(y,x)$, $x=af^i$, and $y=x^\s=a^\s f^i=a^{-1}f^i$.
In particular, $\Ga$ is not arc-transitive, and $S=\{af^i,a^{-1}f^i,(af^i)^{-1},(a^{-1}f^i)^{-1}\}$.

Since $f\in G$ has order $n$,  it follows from Hall's theorem that there exists $b\in W$ such that $h^b\in\l f\r$.
So $f^{b^{-1}}=h^r$ for some $r$. Let $\t=\s^{b^{-1}}$. Then $\t$ centralizes $\l h\r$, and $X=G{:}\l\t\r$.
Moreover, $S^{b^{-1}}=\{ah^{ir},a^{-1}h^{ir},(ah^{ir})^{-1},(a^{-1}h^{ir})^{-1}\}$.
Let $ir\equiv j\,(\mod n)$, and $1\leqslant j\leqslant n-1$. Then $S_j{:}=\{ah^j,a^{-1}h^j,(ah^j)^{-1},(a^{-1}h^j)^{-1}\}$.
Note that $\Ga\cong\Cay(G,S_j)$ is connected.
By \cite[Lemma~6.2]{Hua Zhang}, $(j,n)=1$,
$\Ga_i\cong\Ga_{n-i}$, and if $p^ki\equiv \pm j\,(\mod n)$ for some $k$, then $\Ga_i\cong\Ga_j$.
This completes the proof of  Theorem~\ref{solubles}.
\qed


\begin{thebibliography}{99}
\bibitem{Baik}
Y. G. Baik, Y. Q. Feng, H. S. Sim, M. Y. Xu,
On the normality of Cayley graphs of abelian groups,
{\it Algebra Colloq.} {\bf5} (1998) 297-304.

\bibitem{Magma}
W. Bosma, C. Cannon, C. Playoust, The MAGMA algebra system I: The user language,
{\it J. Symbolic Comput.} {\bf24} (1997) 235-265.



\bibitem{Kuzman}
K. Bo$\check{s}$tjan, Arc-transitive elementary abelian covers of the complete
graph $K_5$, {\it Linear Algebra Appl.} {\bf 433} (2010) 1909-1921.

\bibitem{atlas}
J. H. Conway, R. T. Curtis, S. P. Norton, R. A. Parker, P. Wilson,
Atlas of finite Groups, Oxford University Press, Oxford, 1985.

%\bibitem[Lei-Liu]{Corr}
\bibitem{Corr}
B. P. Corr, C. E. Praeger,
Normal edge-transitive Cayley graphs of Frobenius groups,
{\it J. Algebraic Combin.}  {\bf42} (2015) 803-827.

\bibitem{DF}
A. S. Detinko, D. L. Flannery,
Nilpotent primitive linear groups over finite fields,
{\it Comm. Algebra}. {\bf33} (2005) 497-505.


\bibitem{DM-book}
J. D. Dixon, B. Mortimer,
Permutation Groups, Springer-Verlag, New York, 1996.


\bibitem{Doerk}
K. Doerk, T. Hawkes,
Finite Soluble Groups, Walter de Gruyter Co., Berlin, 1992.

\bibitem{Xin gui}
X. G. Fang, C. H. Li, M. Y. Xu,
On edge-transitive Cayley graphs of valency four,
{\it European J. Combin.} {\bf25} (2004) 1107-1116.


\bibitem{Godsil}
C. D. Godsil, On the full automorphism group of a graph,
{\it Combinatorica}. {\bf1} (1981) 243-256.


\bibitem{Gorenstein}
D. Gorenstein, Finite Groups, Harper and Row, 1968.

\bibitem{L.Kazarin}
L. S. Kazarin, Groups that can be represented as a product of two solvable subgroups,
{\it Comm. Algebra}. {\bf 14} (1986) 1001-1066.

\bibitem{Evgenii I}
E. I. Khukhro, N. Yu. Makarenko, P. Shumyatsky,
Frobenius groups of automorphisms and their fixed points,
{\it Forum Math.} {\bf 26} (2014) 73-112.


\bibitem{Li}
C. H. Li,
The finite vertex-primitive and vertex-biprimitive s-transitive graphs for $s\geqslant4$,
{\it Trans. Amer. Math. Soc.} {\bf 353} (2001) 3511-3529.


\bibitem{flag}
C. H. Li, Finite $s$-arc transitive Cayley graphs and flag-transitive projective planes,
{\it Proc. Amer. Math. Soc.} {\bf 133} (2005) 31-41.

\bibitem{Cai heng}
C. H. Li, Semiregular automorphisms of cubic vertex transitive graphs,
{\it Proc. Amer. Math. Soc.} {\bf 136} (2008) 1905-1910.

\bibitem{Zai Ping Lu}
C. H. Li, Z. Liu, Z. P. Lu,
The edge-transitive tetravalent Cayley graphs of square-free order,
{\it Discrete Math.} {\bf312} (2012) 1952-1967.

\bibitem{Maru}
C. H. Li, Z. P. Lu, D. Maru$\check{s}$i$\check{c}$,
On primitive permutation groups with small suborbits and their orbital graphs,
{\it J. Algebra}. {\bf279} (2004) 749-770.


\bibitem{Hua Zhang}
C. H. Li, Z. P. Lu, H. Zhang,
Tetravalent edge-transitive Cayley graphs with odd number of vertices,
{\it  J. Combin. Theory Ser. B.} {\bf 96} (2006) 164-181.

\bibitem{Xia}
C. H. Li, B. Z. Xia,
Factorizations of almost simple groups with a solvable factor, and Cayley graphs,
arXiv:1408.0350v3~{\bf[math.GR]}.


\bibitem{D.Maru}
D. Maru$\check{s}$i$\check{c}$,
Recent developments in half-transitive graphs,
{\it Discrete Math.} {\bf182} (1998) 219-231.

\bibitem{D.Mar}
D. Maru$\check{s}$i$\check{c}$, R. Nedela,
On the point stabilizers of transitive groups with non-self-paired suborbits of length $2$,
{\it J. Group Theory.} {\bf4} (2001) 19-43.


\bibitem{DD}
D. Maru$\check{s}$i$\check{c}$, R. Nedela,
Maps and half-transitive graphs of valency $4$, {\it European J.Combin.} {\bf 19} (1998) 345-354.

\bibitem{Pan}
J. M. Pan, Y. Liu, Z. H. Huang, C. L. Liu,
Tetravalent edge-transitive graphs of order $p^2q$,
{\it Sci. China Math.} {\bf57(2)}  (2014) 293-302.



\bibitem{p}
P. Poto$\check{c}$nik,
A list of 4-valent 2-arc-transitive graphs and finite faithful amalgams of index (4, 2),
{\it European J. Combin.} {\bf 30} (2009) 1323-1336.


\bibitem{CE}
C. E. Praeger, An O'Nan-Scott theorem for finite quasiprimitive permutation groups and
an application to $2$-arc transitive graphs,
{\it J. London. Math. Soc.} {\bf47} (1992) 227-239.

\bibitem{Praeger}
C. E. Praeger,  Finite normal edge-transitive Cayley graphs,
{\it Bull. Aust. Math. Soc.} {\bf 60(2)} (1999) 207-220.

\bibitem{Suzuki}
M. Suzuki, Group Theory I,
Springer-Verlag, Berlin, New York, 1982.


\bibitem{Xiuyun1}
X. Y. Wang, Y. Q. Feng,
Tetravalent half-edge-transitive graphs and non-normal Cayley graphs,
{\it J. Graph Theory}. {\bf70(2)} (2012) 197-213.


\bibitem{Weiss}
R. M. Weiss, $s$-transitive graphs,
{\it Algebraic Methods in Graph Theory}. {\bf2} (1981) 827-847.


\bibitem{Ming Yao}
M. Y. Xu, Half-transitive graphs of prime-cube order,
{\it J. Algebraic Combin.} {\bf 1(3)} (1992) 275-282.

\bibitem{Xu}
M. Y. Xu, Automorphism groups and isomorphisms of Cayley digraphs,
{\it Discrete Math.} {\bf 182} (1998) 309-319.

\bibitem{Wenqin}
W. Q. Xu, S. F. Du, J. H. Kwak, M. Y. Xu,
$2$-Arc-transitive metacyclic covers of complte graphs,
{\it J. Combin. Theory Ser. B.} {\bf111} (2015) 54-74.

\end{thebibliography}
\end{document}